\documentclass[letterpaper,10pt,conference]{ieeeconf}  %

\IEEEoverridecommandlockouts                              %

\usepackage{balance}

\usepackage{physics}
\usepackage{thmtools}
\usepackage[dvipsnames]{xcolor}
\usepackage{nicefrac}
\usepackage{graphicx,subcaption}
\usepackage{booktabs,multirow,threeparttable}
\usepackage{hyperref}
\usepackage{cleveref}
\usepackage[normalem]{ulem}

\hypersetup{colorlinks = true}

\usepackage{amsmath,amsfonts,bm}

\def\1{\bm{1}}

\def\vw{{\bm{w}}}

\DeclareMathAlphabet{\mathsfit}{\encodingdefault}{\sfdefault}{m}{sl}
\SetMathAlphabet{\mathsfit}{bold}{\encodingdefault}{\sfdefault}{bx}{n}

\newcommand{\R}{\mathbb{R}}

\usepackage{math_macros}
\allowdisplaybreaks

\newtheorem{theorem}{Theorem}
\newtheorem{lemma}[theorem]{Lemma}
\newtheorem{corollary}{Corollary}[theorem]
\newtheorem{proposition}[theorem]{Proposition}

\newtheorem{example}[theorem]{Example}

\newcommand{\splitatcommas}[1]{%
  \begingroup
  \begingroup\lccode`~=`, \lowercase{\endgroup
    \edef~{\mathchar\the\mathcode`, \penalty0 \noexpand\hspace{0pt plus 1em}}%
  }\mathcode`,="8000 #1%
  \endgroup
}

\title{\LARGE \bf
Competitive Control with Delayed Imperfect Information
}

\author{Chenkai Yu\textsuperscript{1}, Guanya Shi\textsuperscript{2}, Soon-Jo Chung\textsuperscript{2}, Yisong Yue\textsuperscript{2}, and Adam Wierman\textsuperscript{2}%
\thanks{\textsuperscript{1}Columbia University, USA. \texttt{cyu26@gsb.columbia.edu}.\endgraf
\textsuperscript{2}Caltech, USA.
\texttt{\{gshi,sjchung,yyue,adamw\}@caltech.edu}.}%
}

\begin{document}

\maketitle
\thispagestyle{empty}
\pagestyle{empty}

\begin{abstract}

This paper studies the impact of imperfect information in online control with adversarial disturbances. In particular, we consider both \emph{delayed} state feedback and \emph{inexact} predictions of future disturbances. We introduce a greedy, myopic policy that yields a constant competitive ratio against the offline optimal policy.  We also analyze the fundamental limits of online control with limited information by showing that our competitive ratio bounds for the greedy, myopic policy in the adversarial setting match (up to lower-order terms) lower bounds in the stochastic setting.

\end{abstract}

\section{INTRODUCTION}

The design and analysis of controllers with imperfect, delayed information is a long-standing challenge for the fields of online learning and robust control.  This paper provides a finite-time analysis of the impact of imperfect information in the context of a disturbed linear system with feedback delay. Specifically, we consider a disturbed online Linear Quadratic Regulator (LQR) optimal control problem with state feedback delay and inexact predictions of future disturbances, governed by $x_{t+1}= Ax_t + B u_t + w_t$, where $x_t$, $u_t$, and $w_t$ are the state, control, and disturbance respectively.  %

A growing literature at the interface of learning and control has emerged in recent years with the goal of designing controllers under various non-asymptotic learning-theoretic criteria, such as regret \cite{dean2018regret,agarwal2019online,simchowitz2020naive,simchowitz2020improper,shi2021meta}, dynamic regret \cite{li2019online,yu2020power,zhang2021regret,lin2021perturbation}, and competitive ratio \cite{shi2020beyond,goel2019online,lin2021perturbation}. However, this line of work has made little progress when it comes to using imperfect prediction information, and has not approached the challenge of delayed feedback at all. %

The task of designing controllers given imperfect predictions has a long history in the control literature, as well as in real-world applications \cite{shi2019neural,lazic2018data}. The basic idea is that, at state $x_t$, one has access to imperfect estimates $\hat{w}_{t+i}$  of the true future disturbances $w_{t+i}$ ($i\geq 0$).  Of particular relevance is robust Model Predictive Control (MPC) or MPC with uncertainty, which focuses on designing MPC policies given inexact future dynamics predictions \cite{robustmpc_survey,cannon2005optimizing}. However, prior work focuses on stability and asymptotic performance, and there are no finite-time performance guarantees known.

The topic of control with delayed feedback is less studied, with some results in the robust control theory literature \cite{zhou1998essentials,kim1999h}, where the focus is on stability and motivation is taken from applications to real-world dynamical systems, e.g., \cite{bejczy1990phantom,shi2021neural}. 
Here, the basic idea is that, at time $t$, one has only observed states up to $x_{t-d}$ ($d\geq 0$).
As in the case of imperfect predictions, there are no finite-time performance guarantees known to the best of our knowledge. 

Thus, while imperfect predictions and delayed feedback are of long-standing importance in both theory and practice, little progress has been made toward obtaining finite-time performance bounds for metrics like regret and competitive ratio under delayed imperfect information due to the technical challenges associated with proving such bounds.

\subsection{Contributions}
We show that a simple, myopic, predictive control has a constant competitive ratio bound in the case of delayed imperfect information (Theorem \ref{thm:delay&inexact_pred}). 
This bound exponentially increases (decreases) in the number of delays (predictions), which highlights the cost associated with delay and the power of predictions, even when they are inexact.
To the best of our knowledge, this result represents the first competitive ratio bounds for either the setting of inexact predictions or delayed feedback. 
We also prove that this bound is tight by showing that
in some systems, the competitive ratio of the optimal online policy can be computed and matches our upper bound.

We would like to emphasize the generality of our result. The model we consider is the general online LQR setting with bounded adversarial disturbance in the dynamics, where only stabilizability is assumed. Further, the prediction errors are assumed to be adversarial. Additionally, our results compare to the globally optimal policies without any constraints (i.e., using the metric competitive ratio), rather than the optimal linear, static policy (i.e., using static regret).

Our result adds further evidence that the structure of LQR allows simple algorithmic ideas to be effective: \cite{simchowitz2020naive} recently proved that the naive exploration is optimal in online LQR adaptive control problem with unknown $\{A,B\}$, and \cite{yu2020power}  proved the classic MPC is near-optimal in online LQR control with exact future predictions.
Combined with the current paper, there is growing evidence that simple, myopic policies that build on MPC are constant-competitive and near-optimal, even in adversarial settings with delayed imperfect information,  which sheds light on key algorithmic principles and fundamental limits in continuous control.

\subsection{Related work} 
There is a growing literature of papers that approach the control of linear dynamical systems with tools and concepts from online learning, with a focus on finite-time non-asymptotic guarantees. This non-asymptotic perspective is important because (1) it can easily integrate with learning theory and (2) finite-time optimality analysis is crucial for many real-world systems with finite horizons. Within this literature, most work focuses on the design of controllers with low regret \cite{dean2018regret,agarwal2019online,simchowitz2020naive,yu2020power,zhang2021regret,nonhoff2020online}.
The study of competitive ratio has received increased attention recently~\cite{goel2019online,shi2020beyond,goel2021competitive,lin2021perturbation}. Again, note that competitive ratio compares with the offline optimal policy without any constraints while regret compares with the optimal static policy in a specific policy class. The most general such result to this point is~\cite{lin2021perturbation} which provides competitive bounds of a predictive controller in a time-varying setting. However, the bounds in \cite{lin2021perturbation} require a sufficiently large amount of \emph{exact} predictions.

In these lines of work, very few papers focus on settings where the controller has access to predictions of future disturbances. Among these papers, \cite{yu2020power,lin2021perturbation} focus on settings where predictions are exact. To the best of our knowledge, the only result for inexact predictions is studied in \cite{zhang2021regret}, however \cite{zhang2021regret} only gives dynamic regret results which are data-dependent. In this paper we provide stronger, data-independent results for competitive ratio.
Inexact dynamic predictions have received attention in the robust MPC community (e.g., \cite{robustmpc_survey,cannon2005optimizing}), focusing on stability and asymptotic performance analysis. In contrast, this paper focuses on non-asymptotic analysis from a learning theory perspective. 
Even outside of control in the related area of online optimization, when predictions are considered, they are typically assumed to be exact \cite{lin2019online}. One exception is \cite{chen2015online}, which uses a less general model of prediction error than the current paper, and the connection to control is unclear. 

In contrast to the literature on predictions, there is no work providing non-asymptotic guarantees (either regret or competitive ratio) of policies subject to delayed feedback. The issue has received considerable attention in the robust control community \cite{kim1999h,bejczy1990phantom}, but the focus is typically on closed-loop stability and no finite-time performance bounds exist to the best of our knowledge.

\section{MODEL} \label{section:model}
We consider an online \emph{Linear Quadratic Regulator (LQR)} optimal control problem with adversarial disturbances in the dynamics. In particular, we consider a linear system initialized with $x_0 \in \R^n$ and controlled by $u_t \in \R^m$ at each step $t \in \{0, 1, \dots, T-1\}$ where $T$ is the total length of the problem. The system dynamics is governed by:
\[x_{t+1} = A x_t + B u_t + w_t,\]
where $w_t$ is the disturbance. We assume that $w_t$ is bounded, i.e., $\norm{w_t} \le r$. 
The goal is to minimize the following cost:
\begin{equation} \label{eq:cost_function}
  J = \sum_{t=0}^{T-1} (x_t^\trp Q x_t + u_t^\trp R u_t) + x_T^\trp \Qf x_T,
\end{equation}
given matrices $A, B, Q, R, \Qf$. We consider an \emph{online} setting where an adversary \emph{adaptively} selects $\{w_t\}_{t=0}^{T-1}$, and the controller (also adaptively) makes the decision $u_t$ at every time step $t$, potentially based on \emph{delayed imperfect information} %
(discussed later in this section).

We make the standard assumptions that $Q, \Qf \succeq 0$, $R \succ 0$, and $(A, B)$ is stabilizable \cite{goel2020power,yu2020power}, i.e., $\exists K_0 \in \R^{m \times n}$ such that $\rho(A - B K_0) < 1$. We further assume $(A^\trp, Q)$ is stabilizable to guarantee the stability of the closed-loop \cite{anderson2007optimal}. This assumption is more general than the standard assumption that $Q \succ 0$ because $Q \succ 0$ implies the stabilizability of $(A^\trp, Q)$. 

Throughout this paper, we use $\rho(\cdot)$ to denote the spectral radius of a matrix and $\norm{\cdot}$ to denote the 2-norm of a vector or the spectral norm of a matrix.

Note that many important problems can be considered as special cases of the above model. One motivating example is the Linear Quadratic (LQ) tracking problem \cite{anderson2007optimal}.

\begin{example}[Linear quadratic tracking]
The LQ tracking problem is defined via dynamics $x_{t+1}=Ax_t+Bu_t+\Tilde{w}_t,$ and cost function
$J=\sum_{t=0}^{T-1} (x_{t+1}-y_{t+1})^\trp Q (x_{t+1}-y_{t+1}) + u_t^\trp R u_t,$
where $\{y_t\}_{t=1}^T$ is the desired trajectory to track. 
\end{example}

To fit LQ tracking into our model, let $\Tilde{x}_t=x_t-y_t$. Then, we get
$J=\sum_{t=0}^{T-1} \Tilde{x}_{t+1}^\trp Q \Tilde{x}_{t+1} + u_t^\trp R u_t$
and $\Tilde{x}_{t+1} = A \Tilde{x}_t + B u_t + w_t$,
which is a LQR control problem with disturbance $w_t = \tilde{w}_t + Ay_t - y_{t+1}$ in the dynamics.  Note that in many LQ tracking problems, delayed observations and imperfect predictions are fundamental challenges \cite{anderson2007optimal}.

\subsection{Delayed Imperfect Information}
\label{section:pred&delay}

The LQR optimal control problem introduced above is typically studied without \emph{predictions} or \emph{delays}.
In the classic setting, at each time $t = 0, 1, 2, \dots$, the controller observes $x_t$ and then decides $u_t$ without knowing $w_t$.
Thus, $u_t$ is a function of all previous information: $u_t = \pi_t(\splitatcommas{x_0, x_1, \dots, x_{t}, u_0, u_1, \dots, u_{t-1}})$.
In an equivalent formulation, the controller is given $x_0$ to start and, at each time $t$, obtains the previous disturbance $w_{t-1}$ (if $t \ge 1$) and then decides $u_t$.
Thus, $u_t = \pi_t(\splitatcommas{x_0, w_0, w_1, \dots, w_{t-1}, u_0, u_1, \dots, u_{t-1}})$.

The motivation for this paper is that 
in many real-world problems (e.g., \cite{shi2019neural,lazic2018data}), predictions of future information are available and using them is crucially important, even though they are typically noisy and delayed. For example, data-driven model-based control is a prominent and successful approach where it is crucial to consider model mismatch due to statistical learning error. Moreover, in many situations, there is state feedback delay in the system, so the controller has to make the decision $u_t$ before the current state $x_t$ is observed. The existence of both imperfect prediction and feedback delay leads to considerable difficulty.

Formally, we model the delayed inexact predictions %
as follows.
At each step $t$, the revealed information is:
\[x_0,\ u_0, \dots, u_{t-1},\ w_0, \dots, w_{t-d-1},\ \hat w_{t-d|t}, \dots, \hat w_{T-1|t},\]
or equivalently,
\[x_0,\ u_0, \dots, u_{t-1},\ x_1, \dots, x_{t-d},\ \hat w_{t-d|t}, \dots, \hat w_{T-1|t},\]
where $d \ge 0$ is the length of feedback delay, and $\hat w_{s|t}$ is the %
prediction of $w_s$ at time $t$. %

We define $e_{s|t} = w_{s} - \hat w_{s|t}$ as the estimation error, and we assume that the predictor satisfies $\norm{e_{t-d+i|t}} \le \epsilon_i \norm{w_{t-d+i}}$ for all $i \ge 0$ and $t \ge 0$, where $\epsilon_i$ is a given parameter that measures the prediction quality at $i$ steps into the unknown.
Although predictions are available for every time step, those for far future may have bad quality, i.e., $\epsilon_i$ is typically large for large $i$. Therefore, a good control policy may not use all predictions in the same way: using only the predictions with smaller estimation error may yield better performance. For example, if $\epsilon_i > 1$, we can simply let $\hat w_{t-d+i|t} = 0$, which is a better prediction with $\epsilon_i = 1$. In general, the error bounds have to be multiplicative rather than additive, if one pursues a competitive ratio guarantee.
Consider the case where all disturbances $w_t$ are zero, but there are nonzero prediction errors. The optimal offline policy incurs zero cost, while any online policy that uses the (wrong) predictions incurs nonzero cost, hence leading to an infinite ratio.

Our setting of \emph{delayed imperfect information} generalizes many existing settings in the study of LQR control. The classic setting is the special case where $d = 0$ and $\hat w_{t+i|t} = 0$ for all $i$ and $t$. The offline optimal setting considered by \cite{goel2020power,foster2020logarithmic} is the special case where $d = 0$ and $\epsilon_i = 0$ for all $i$.
The setting considered by \cite{yu2020power}, where $k$ \emph{exact} predictions without delay are available, corresponds to $d = 0$, $\epsilon_i = 0$ for $0 \le i \le k-1$ and $\hat w_{t+i|t} = 0$ for all $t$ and $i \ge k$.

\subsection{Competitive Ratio}
We measure the performance using the \emph{competitive ratio},
which bounds the worst-case ratio of the cost of an online policy ($\Alg$) to the cost of the optimal offline policy ($\Opt$) with perfect knowledge of $\{w_t\}_{t=0}^{T-1}$. Formally, we study the so-called \emph{weak} competitive ratio, which allows for an additive horizon-independent constant factor. We say that a policy is (weakly) $c$-competitive if, given $A$, $B$, $Q$, $R$, $\Qf$ and $r$, for any adversarially and adaptively chosen disturbances $\{w_t\}_{t=0}^{T-1}$, we have $\Alg \le c \, \Opt + \kappa$, where $\kappa$ is a constant, independent of $T$. We say that the algorithm is \emph{constant competitive} when $c$ is a constant independent of $T$.

While there has been considerable success in recent years designing control policies that are no-regret,
e.g., \cite{dean2018regret,agarwal2019online,yu2020power}, there have been very few examples of constant competitive controllers. The few results that exist, e.g., \cite{shi2020beyond,goel2019online}, tend to have more restrictive assumptions on the dynamics and/or disturbances.
Existing results on regret typically focus on \emph{static} regret, which compares a policy to the optimal offline \emph{static} policy
in a specific policy class (e.g., linear policies~\cite{agarwal2019online}),
while the \emph{dynamic} regret and the competitive ratio compare a policy to the general optimal offline policy that may potentially be non-linear and non-static.
Although it is possible for an online policy to have sublinear \emph{static} regret, in general, the optimal static linear policy may have cost that is an order-of-magnitude larger than the optimal offline cost
\cite{shi2020beyond}.
On the other hand, it has been shown \cite{yu2020power} that the minimum \emph{dynamic} regret of an online policy can be linear in $T$ even with random disturbances. Thus, it is natural to consider the competitive ratio and pin down the multiplicative constant incurred from not knowing the future.

\section{A MYOPIC PREDICTIVE POLICY}
\label{sec:myopic_policy}
In this paper, we study a myopic predictive control policy that is a natural and myopic variant of model predictive control (MPC). When there are predictions, but no delays, MPC is a popular and successful approach \cite{lazic2018data,camacho2013model}. In fact, \cite{yu2020power} recently showed that MPC has a near-optimal dynamic regret in the case of exact predictions and no delay. 
Note that this paper aims to show simple/standard algorithms (either standard MPC or its natural extension) could be competitive and the competitive ratio bounds could be tight, rather than propose a new algorithm.

Suppose the controller uses $k$ predictions. At each time $t$, the controller optimizes based on $\splitatcommas{x_t, \hat w_{t|t}, \dots, \hat w_{t+k-1|t}}$:
\begin{align}
&\begin{aligned}
(u_t, & \dots, u_{t+k-1}) = \\
&\arg\min_u \bigg( \sum_{i=t}^{t+k-1} (x_i^\trp Q x_i + u_i^\trp R u_i) \label{eq:mpc_step} 
 + x_{t+k}^\trp \tilde\Qf x_{t+k} \bigg), 
\end{aligned} \notag \\
&\text{s.t. } x_{i+1} = A x_i + B u_i + \hat w_{i|t}, \ \forall i = t, \dots, t+k-1.
\end{align}
This optimization is myopic in the sense that it assumes that the length of the problem is $k$ instead of $T$ and only uses predicted future disturbances within those $k$ steps. The terminal cost matrix $\tilde\Qf$ in \eqref{eq:mpc_step} may or may not be the same as the terminal cost matrix $\Qf$ of the original problem \eqref{eq:cost_function}, and can be viewed as a hyperparameter. Similarly, $k$ is also a hyperparameter. Larger $k$ is not necessarily better because the predictions in the far future may have very large errors. In this paper, we let $\tilde\Qf = P$, where $P$ is the solution of the discrete algebraic Riccati equation (DARE):
\begin{equation} \label{eq:DARE}
  P = Q + A^\trp P A - A^\trp P B (R + B^\trp P B)^{-1} B^\trp P A.
\end{equation}
The output of \eqref{eq:mpc_step} is $k$ control actions corresponding to time $t, t + 1, \dots, t + k - 1$, respectively, but only the first ($u_t$) is applied to the system. The rest (i.e., $u_{t+1}, \dots, u_{t+k-1}$) are discarded.  The explicit solution of \eqref{eq:mpc_step} is given below \cite{yu2020power}.
\begin{proposition} \label{prop:mpc_policy}
  The policy defined by \eqref{eq:mpc_step} at time $t$ is:
  \[u_t = -(R + B^\trp P B)^{-1} B^\trp \left(P A x_t + \sum_{i=0}^{k-1} {F^\trp}^i P \hat w_{t+i|t} \right),\]
  where $F = A - B (R + B^\trp P B)^{-1} B^\trp P A =: A - BK$.
\end{proposition}

The closed loop is stable since $\rho(F) < 1$ \cite{yu2020power}. As stated above, the policy does not directly apply to the case of delayed imperfect information. To adapt it, we consider two cases: (i) when the number of predictions available is longer than the feedback delay, i.e., $k \ge d$, and (ii) when the delay is longer than the number of predictions available, i.e., $k < d$.

When $k \ge d$, the extension is perhaps straightforward. Here, although the controller does not know the current state $x_t$, it knows $x_{t-d}$ and $\hat w_{t-d|t}, \dots, \hat w_{t-1|t}$.  Thus, it can estimate the current state. This means that it is possible to simply use this estimation, $\hat x_{t|t}$, as a replacement for $x_t$ in the algorithm, which yields the following:
\begin{equation} \label{eq:myopic_policy_k_gt_d}
  u_t = -(R + B^\trp P B)^{-1} B^\trp \bigg(P A \hat x_{t|t} 
  + \sum_{i=0}^{k-d-1} {F^\trp}^i P \hat w_{t+i|t} \bigg).
\end{equation}

When $k < d$, the extension is not as obvious. In this setting, the quality of the predictions is poor enough that it is better not to use the predictions to estimate the current state. Thus, one cannot simply estimate the current state and run classic MPC. In this case, the key is to view (classic) MPC from a different perspective: MPC locally solves an optimal control problem by treating known disturbances (predictions) as exact, and treating unknown disturbances as zero \cite{yu2020power,camacho2013model,zhang2021regret}.
Following this philosophy, in the case when predictions are not enough to be used to estimate the current state, we can instead assume that unknown disturbances are exactly zero. The following theorem derives the optimal policy under this ``optimistic'' assumption.

\begin{theorem}
\label{thm:myopic_policy_delay}
  Suppose there are $d$ delays and $k$ exact predictions with $k < d$. Assume all used predictions are exact and other disturbances (with unused predictions) are zero. The optimal policy at time $t$ is:
  \begin{multline} \label{eq:myopic_policy_k_lt_d}
    u_t = -(R + B^\trp P B)^{-1} B^\trp P A \bigg(A^{d-k} \hat x_{t-d+k|t} \\
    + \sum_{i=0}^{d-k-1} A^i B u_{t-1-i}\bigg).
  \end{multline}
\end{theorem}

In other words, the policy in \eqref{eq:myopic_policy_k_lt_d} first obtains the greedy estimation $\hat x_{t-d+k|t}$ using predictions $\hat w_{t-d|t},\dots,\hat w_{t-d+k-1|t}$, and then estimates the current state by treating $w_{t-d+k}=\dots=w_{t-1}=0$.
In fact, instead of treating them as zero, we can impose other values or distributions on those disturbances. This would generalize \Cref{thm:myopic_policy_delay} to a broader class of policies.

To summarize the two cases above, the myopic generalization of MPC we study in this paper is described as follows.
Suppose we want to use $k$ predictions. If $k \ge d$, then we estimate the current state $x_t$ and apply \eqref{eq:myopic_policy_k_gt_d}. If $k < d$, then we estimate the state at time $t - d + k$ and apply \eqref{eq:myopic_policy_k_lt_d}. If fact, the two cases coincide when $k = d$.

\section{PERFORMANCE BOUNDS}

Our main result provides bounds on the competitive ratio for the policy defined in \Cref{sec:myopic_policy} in the case of inexact delayed predictions.   We present our general result below and then discuss the special cases of (i) exact predictions and no delay, (ii) inexact predictions and no delay, and (iii) delay but no access to predictions.  The special cases illustrate the contrast between inexact and exact predictions as well as the impact of delay.

\begin{theorem}[Main result] \label{thm:delay&inexact_pred}
  Let $c = \norm{P} \|P^{-1}\| (1 + \norm{F})$ and $H = B (R + B^\trp P B)^{-1} B^\trp$. Suppose there are $d$ steps of delays and the controller uses $k$ predictions. When $k \ge d$,
  \[\begin{split}
    \Alg & \le \bigg[\frac{\left(c \sum\limits_{i=0}^{d-1} \epsilon_i \|A^{d-i}\| + c \sum\limits_{i=d}^{k-1} \epsilon_i \|F^{i-d}\| + \|F^{k-d}\|\right)^2}{\norm{H}^{-1} \lambda_{\min}(P^{-1} - F P^{-1} F^\trp - H)} \\
    & \hspace{5cm} + 1 \bigg] \Opt + O(1).
  \end{split}\]
  When $k \le d$,
  \[\begin{split}
    \Alg & \le \bigg[\frac{\bigg(c \sum\limits_{i=0}^{k-1} \epsilon_i \|A^{d-i}\| + c \sum\limits_{i=k}^{d-1} \|A^{d-i}\| + 1 \bigg)^2}{\norm{H}^{-1} \lambda_{\min}(P^{-1} - F P^{-1} F^\trp - H)} \\
    & \hspace{5cm} + 1 \bigg] \Opt + O(1).
  \end{split}\]
  The $O(1)$ is with respect to $T$. It may depend on the system parameters $A$, $B$, $Q$, $R$, $\Qf$ and the range of disturbances $r$, but not on $T$. When $\Qf = P$ the $O(1)$ is zero.
\end{theorem}

The two cases in the theorem correspond to the two cases in the algorithm: when predictions are of high enough quality to allow estimation of the current state and when they are not. Note that the closed-loop dynamics is stable, i.e., $\rho(F)=\rho(A-BK)<1$. Therefore, there exists a constant $\gamma$ such that $\norm{F^i} \le \gamma(\frac{\rho(F)+1}{2})^i$ for all $i \ge 1$ from Gelfand’s formula.

In the first case, we see that the quality of predictions in the near future has more impact, especially when $\rho(A)>1$. In the second case, we see that the amount of delay $d$ exponentially increases the bound if $\epsilon_i > 0$ and $\rho(A)>1$.
We explore further insights by looking at special cases in the following subsections.
Before moving on, we provide an overview of the proof of Theorem \ref{thm:delay&inexact_pred}.

\subsection{Proof Sketch for Theorem \ref{thm:delay&inexact_pred}}

We first prove \Cref{thm:delay&inexact_pred} in the case $\Qf = P$. In this case, the $O(1)$ is not needed.  Then, we analyze the impact of the terminal cost $\Qf$ and show that it introduces at most an $O(1)$ additional cost. The full proof can be found in the appendix.

\begin{lemma} \label{lemma:when_Qf=P}
  The conclusion of \Cref{thm:delay&inexact_pred} holds if $\Qf = P$.
\end{lemma}

The proof of the result in the case of $\Qf = P$ follows from a novel difference analysis of the quadratic cost-to-go functions.  Here, we focus on the second part of the proof, i.e., reducing the case of $\Qf \ne P$ to the case of $\Qf = P$. To that end, let $\Alg(X)$ be the cost of our algorithm when the terminal cost is $X$ and, similarly, let $\Opt^Y(X)$ be the cost of the policy that is optimal for terminal cost $Y$ when the terminal cost is actually $X$.

Our analysis proceeds by first bounding the impact of the terminal cost on the gap between the algorithm and the optimal cost via the following lemma.

\begin{lemma} \label{lemma:cost_gap_O(1)}
  For any algorithm,
  \[\Alg(\Qf) - \Opt^P(\Qf) \le \Alg(P) - \Opt^P(P) + O(1).\]
\end{lemma}

Then, we prove that the terminal cost only has an $O(1)$ impact on the optimal cost using the following lemma.

\begin{lemma} \label{lemma:terminal_cost_Opt}
  The followings are equal up to $O(1)$ difference: $\Opt^P(P),\ \Opt^P(\Qf),\ \Opt^\Qf(\Qf)$.
\end{lemma}

Together, these two lemmas imply that
\begin{align*}
  & \Alg(\Qf) - \Opt^\Qf(\Qf) \\
  & \le \Alg(\Qf) - \Opt^P(\Qf) + O(1) \\
  & \le \Alg(P) - \Opt^P(P) + O(1) \\
  & \le \frac{\Alg(P) - \Opt^P(P)}{\Opt^P(P)} \Opt^P(P) + O(1) \\
  & \le \frac{\Alg(P) - \Opt^P(P)}{\Opt^P(P)} \Opt^\Qf(\Qf) + O(1).
\end{align*}
Thus, we can complete the proof by concluding that
\[\Alg(\Qf) \le \frac{\Alg(P)}{\Opt^P(P)} \Opt^\Qf(\Qf) + O(1).\]

\subsection{Exact Predictions Without Delay}
In the sections that follow, we explore special cases of Theorem \ref{thm:delay&inexact_pred} in order to highlight the impact of inexact predictions and delay.  First, we present the special case of $k$ accurate, exact predictions and no feedback delays. Formally, we have $d = 0$, $\epsilon_i = 0$ for $0 \le i \le k - 1$ and $\hat w_{t+i|t} = 0$ for all $t$ and $i \ge k$.

The main result for this setting is given below. It directly follows from the $k \ge d$ case of \Cref{thm:delay&inexact_pred}. 

\begin{theorem} \label{thm:pred_CR}
  Suppose there are $k$ exact predictions and no feedback delay.  Then:
  \[\Alg \le \left[1 + \frac{\|F^k\|^2 \|H\|}{\lambda_{\min}(P^{-1} - F P^{-1} F^\trp - H)}\right] \Opt + O(1).\]
\end{theorem}

\begin{figure*}[htb]
  \centering
  \subcaptionbox{$k = 0$.}{\includegraphics[height = 0.25 \hsize]{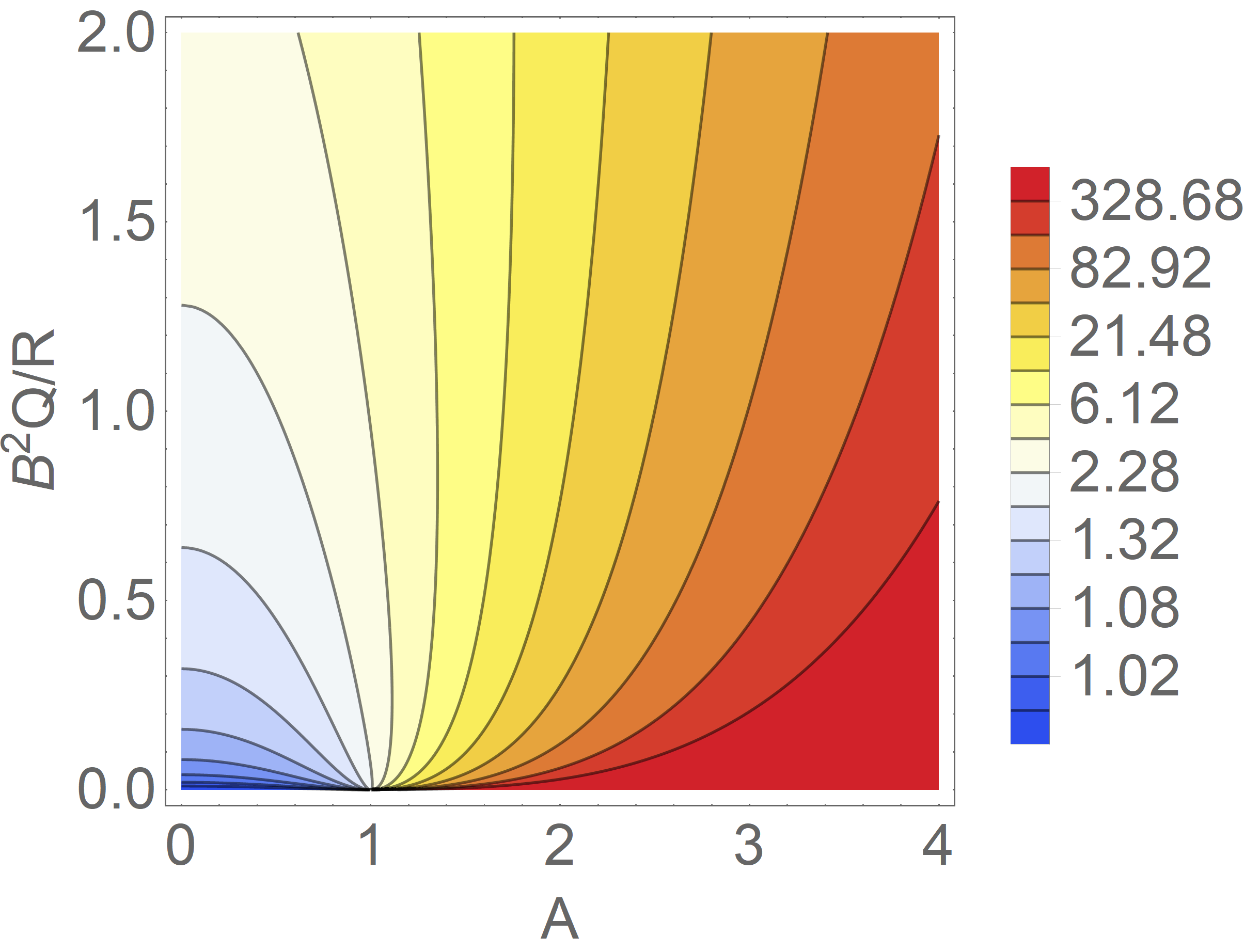}}
  \hfill
  \subcaptionbox{$k = 1$.}{\includegraphics[height = 0.25 \hsize]{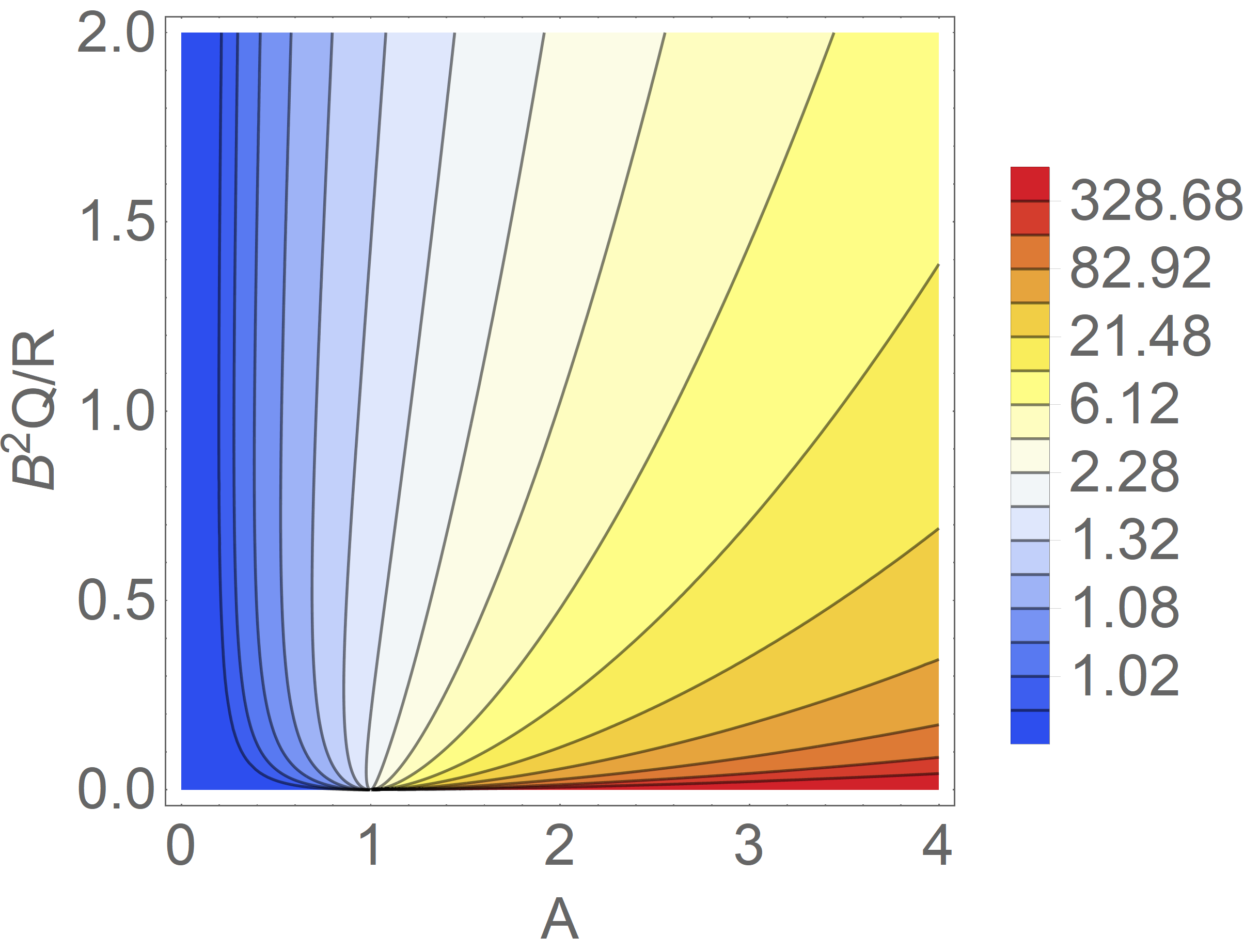}}
  \hfill
  \subcaptionbox{$k = 3$.}{\includegraphics[height = 0.25 \hsize]{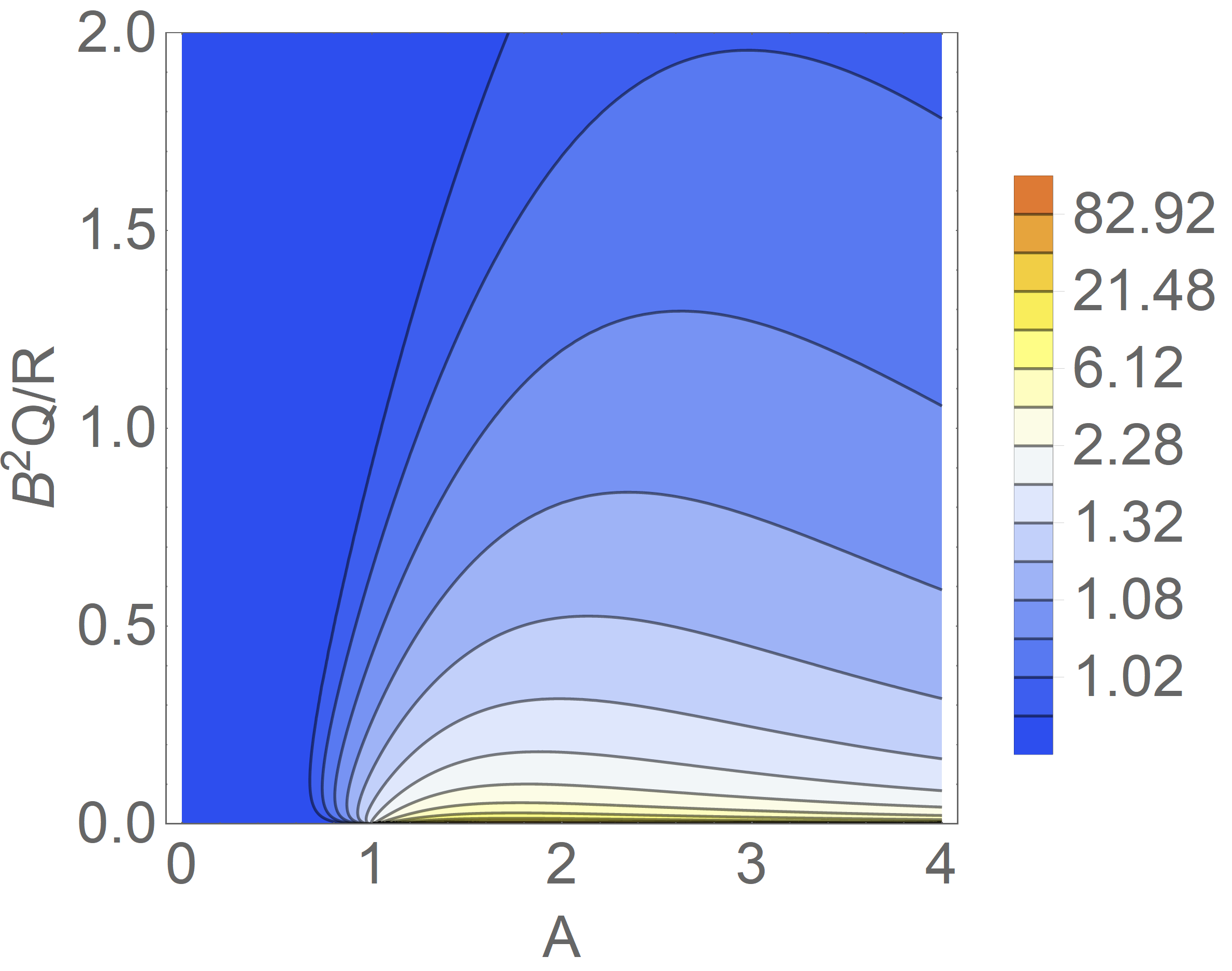}}
  \caption{Illustration of the competitive ratio bound in \Cref{thm:1-d_pred_CR}.
  The system is one-dimensional ($n = m = 1$) without delays ($d = 0$). With no predictions ($k = 0$), the bound is small only if both $A$ and $B^2 Q / R$ are small. When $k = 1$, it is small if $A$ is small or $B^2 Q / R$ is large. When $k = 3$, it is small if $A$ is either small or large, or if $B^2 Q / R$ is large.
  The bound is tight in the sense that in some systems, it equals the relative cost of the optimal online policy.
  }
  \label{fig:CR_bound_1-dim}
\end{figure*}

Thus, the competitive ratio exponentially decreases as $k$ goes up. To illustrate how the primary parameters $A$, $B$, $Q$ and $R$ affect the competitive ratio in this case, it is useful to consider the case when $n=m=1$, as shown in both \Cref{thm:1-d_pred_CR} and \Cref{fig:CR_bound_1-dim}.

\begin{corollary} \label{thm:1-d_pred_CR}
  Assume there are $k$ exact predictions and no feedback delay, and let $n = m = 1$ and $\Qf = P$. Then,
  \begin{alignat*}{2}
    \frac{\Alg}{\Opt} & \le 1 + \frac{2 A^{4-2k}}{B^2 Q / R} \quad & & \text{if } A^2 \gg B^2 Q / R + 1, \\
    \frac{\Alg}{\Opt} & \le 1 + \frac{A^{2k}}{(B^2 Q / R)^{2k-1}} \quad & & \text{if } B^2 Q / R \gg A^2 + 1.
  \end{alignat*}
\end{corollary}

Interestingly, in this case, the competitive bound only depends on $A^2$ and $B^2 Q / R$. It does not depend on the sign of $A$, nor on $B$, $Q$ or $R$ as long as $B^2 Q / R$ is fixed.  Further, when $k \ge 3$, we see that the competitive ratio is small if $B, Q$ are small, $R$ is large, or $A$ is either very large or very small. However, when $k = 0$ or $1$, a large $A$ can result in a large competitive ratio. When $k = 0$, a large value of $B^2 Q / R$ also results in a large competitive ratio. We see below that this is similar to the case of delay (see \Cref{section:delay_special_case}).

\begin{proposition}
  \Cref{thm:pred_CR} is tight in the sense that there exist systems where the competitive ratio of the optimal online algorithm is $1 + \Theta(\|F^k\|^2)$.
\end{proposition}

\subsection{Inexact Predictions Without Delay}
We next consider the case where predictions are inexact, but there is no feedback delay. The contrast with the previous section highlights the impact of prediction error.

As discussed in \Cref{section:pred&delay}, the controller should optimize $k$ to utilize predictions with smaller estimation errors while avoiding the use of those with larger errors. The following directly follows from the $d=0$ case of \Cref{thm:delay&inexact_pred} and reduces to the exact case when $\epsilon_i = 0$. %

\begin{theorem} \label{thm:inexact_pred_CR}
  Suppose there are $k$ inexact predictions and no feedback delay. Then,
  \[\Alg \le \bigg[ \frac{\norm{H} \left(c \sum_{i=0}^{k-1} \epsilon_i \|F^i\| + \|F^k\|\right)^2 }{\lambda_{\min}(P^{-1} - F P^{-1} F^\trp - H)} + 1 \bigg] \Opt + O(1).\]
\end{theorem}

This subsection differs from the previous one in that the controller can minimize the bound in \Cref{thm:inexact_pred_CR} with respect to $k$. We characterize this optimization in the following result in 1-d systems, and also provide simulation evidence in \Cref{section:simulation} (see \Cref{fig:inexact}).

\begin{corollary}
  Suppose there are $k$ inexact predictions and no feedback delay. Assume $n = m = 1$. Given non-decreasing $\{\epsilon_i\}$, to minimize the competitive ratio bound in \Cref{thm:inexact_pred_CR}, the optimal number $k$ of predictions to use is such that:
  \[\epsilon_{k-1} < \frac{1 - \abs{F}}{1 + \abs{F}} < \epsilon_k.\]
\end{corollary}

The following 1-d setting highlights the dependence of the competitive ratio on the system parameters. 

\begin{corollary} \label{thm:1-d_inexact_CR}
  Assume there are $k$ inexact predictions and no feedback delay, and let $n = m = 1$ and $\Qf = P$.
  
  If $A^2 \gg B^2 Q / R + 1$,
  \[\frac{\Alg}{\Opt} \le 1 + \frac{2 A^4}{B^2 Q / R} \left(\sum_{i=0}^{k-1} \frac{\epsilon_i}{\abs{A}^i} + \frac{1}{\abs{A}^k}\right)^2.\]
  If $B^2 Q / R \gg A^2 + 1$,
  \[\frac{\Alg}{\Opt} \le 1 + \frac{B^2 Q}{R} \left(\sum_{i=0}^{k-1} \frac{\epsilon_i \abs{A}^i}{(B^2 Q / R)^i} + \frac{\abs{A}^k}{(B^2 Q / R)^k}\right)^2.\]
\end{corollary}

The dependence of the competitive ratio on $A$, $B$, $Q$, $R$ is similar to the case of exact predictions. In particular, we find that the prediction quality in the near future is (exponentially) more important than further in the future, which is consistent with the robust MPC literature \cite{cannon2005optimizing}.

In the exact prediction case we show that \Cref{thm:pred_CR} is tight with respect to $\norm{F^k}$. In contrast, in the inexact case the tightness of $\epsilon_i$ and $\norm{F^i}$ in \Cref{thm:inexact_pred_CR} remains as an open question.

\subsection{Delay Without Predictions} \label{section:delay_special_case}

The last special case we consider is the case with delays but no (usable) predictions.  This case separates the impact of delay from that of predictions.  Here, $\hat w_{t-d+i|t} = 0$ for all $t$ and $i \ge 0$.
When $k \le d$, via \Cref{thm:delay&inexact_pred} we have that:
\[\Alg \le \bigg[ \frac{\norm{H} \left(c \sum_{i=1}^d \|A^i\| + 1 \right)^2}{\lambda_{\min}(P^{-1} - F P^{-1} F^\trp - H)} + 1 \bigg] \Opt + O(1).\]
Depending on whether the spectral radius $\rho(A) < 1$, there are two simplifications one can make: (i) if $\rho(A) < 1$, then $\|A^i\| \le \kappa a^i$ for $a = (\rho(A) + 1) / 2 < 1$ for a constant $\kappa$, and (ii) if $\rho(A) > 1$, $\|A^i\| \le \|A\|^i$.

\begin{theorem} \label{thm:delay_CR}
  Suppose there are $d$ delays and no predictions are available.  If $\rho(A) < 1$, then the competitive ratio is bounded by a constant irrelevant to the length of delay:
  \[\begin{split}
    \Alg
    & \le \bigg[ \frac{\norm{H} \Big(c \kappa \dfrac{a}{1 - a} + 1 \Big)^2}{\lambda_{\min}(P^{-1} - F P^{-1} F^\trp - H)} + 1 \bigg] \Opt + O(1).
  \end{split}\]
  If $\rho(A) > 1$, then the competitive ratio bound grows exponentially fast with the of number of delay steps:
  \[\Alg \le \bigg[ \frac{\norm{H} \bigg(c \dfrac{\norm{A}^{d+1} - \norm{A}}{\norm{A} - 1} + 1 \bigg)^2}{\lambda_{\min}(P^{-1} - F P^{-1} F^\trp - H)} + 1 \bigg] \Opt + O(1).\]
\end{theorem}

As in the previous subsections, it is useful to consider the one-dimensional case to get insights about the impact of the system parameters.

\begin{corollary} \label{thm:1-d_delay_CR}
  Assume there are $d$ delays and no predictions. Let $n = m = 1$ and $\Qf = P$. Then,
  \begin{alignat*}{2}
    \frac{\Alg}{\Opt} & \le 1 + \frac{2 A^{4+2d}}{B^2 Q / R} \ & & \text{if } A^2 \gg B^2 Q / R + 1, \\
    \frac{\Alg}{\Opt} & \le 1 + \frac{A^{2d+2} B^2 Q / R}{(\abs{A} - 1)^2} \ & &\text{if } B^2 Q / R \gg A^2 + 1, \, \abs{A} > 1, \\
    \frac{\Alg}{\Opt} & \le 1 + \frac{B^2 Q / R}{(1 - \abs{A})^2} \ & & \text{if } \abs{A} < 1.
  \end{alignat*}
\end{corollary}

Contrary to the case with $k \ge 3$ inexact predictions, when there are delays, a large value of $B^2 Q / R$ or $A$ does not lead to a small competitive ratio. Instead, in the case of feedback delay, it results in a large competitive ratio. This is consistent with results from robust control theory: the less stable the open loop is ($|A|$ is larger), the more impact delay has \cite{zhou1998essentials}.

\begin{proposition}
  \Cref{thm:delay_CR} is tight in the sense that there exist systems such that the competitive ratio of the optimal online algorithm is at least $1+\Theta(\|A^d\|^2)$.
\end{proposition}

\section{NUMERICAL EXAMPLES} \label{section:simulation}

\begin{figure*}[t]
  \centering
  \subcaptionbox{$d = 5,\ k = 0$.}{\includegraphics[width = 0.18 \hsize]{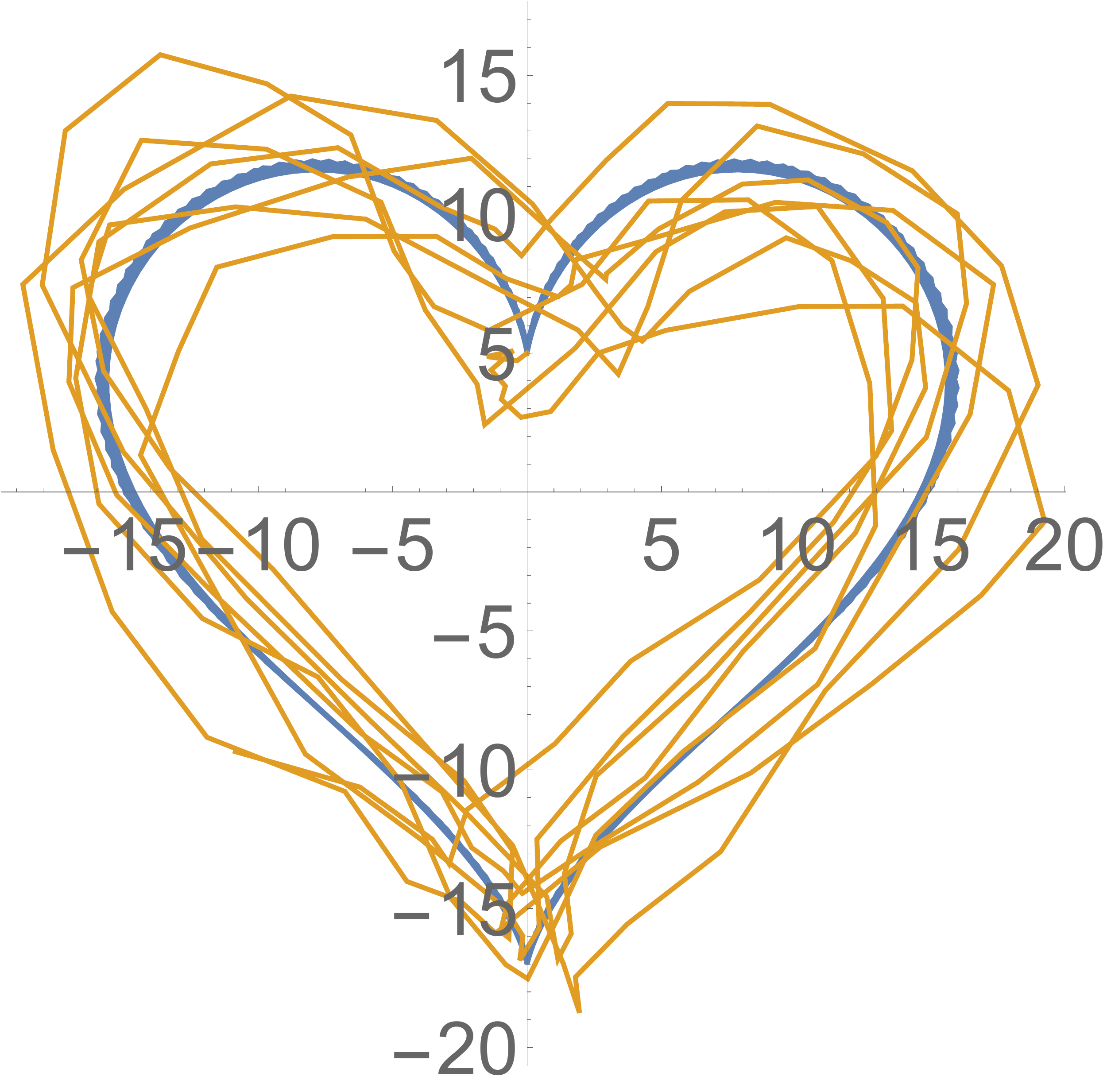}}
  \hfill
  \subcaptionbox{$d = 1,\ k = 0$.}{\includegraphics[width = 0.18 \hsize]{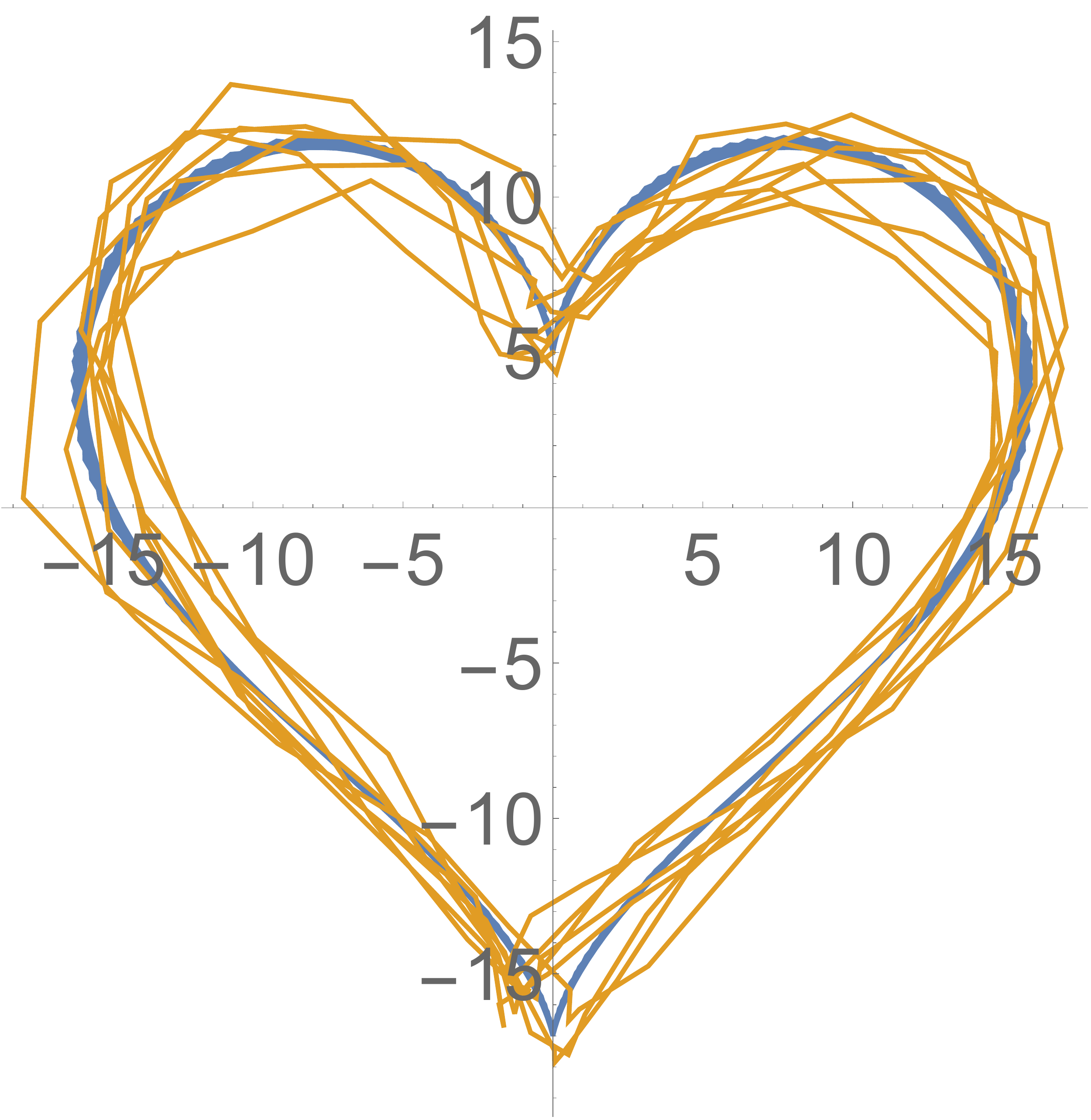}}
  \hfill
  \subcaptionbox{$d = 0,\ k = 1$.}{\includegraphics[width = 0.18 \hsize]{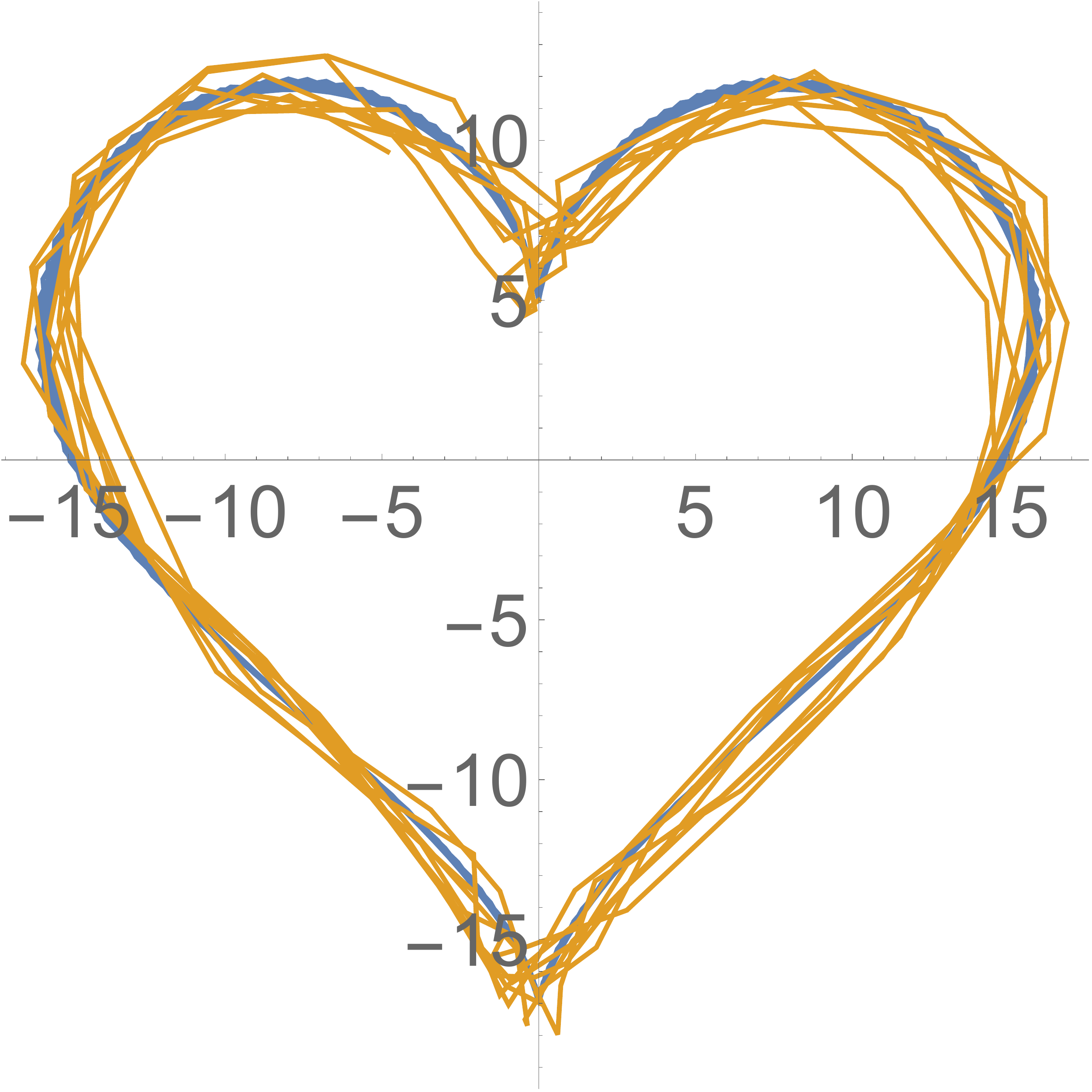}}
  \hfill
  \subcaptionbox{$d = 0,\ k = 8$.}{\includegraphics[width = 0.18 \hsize]{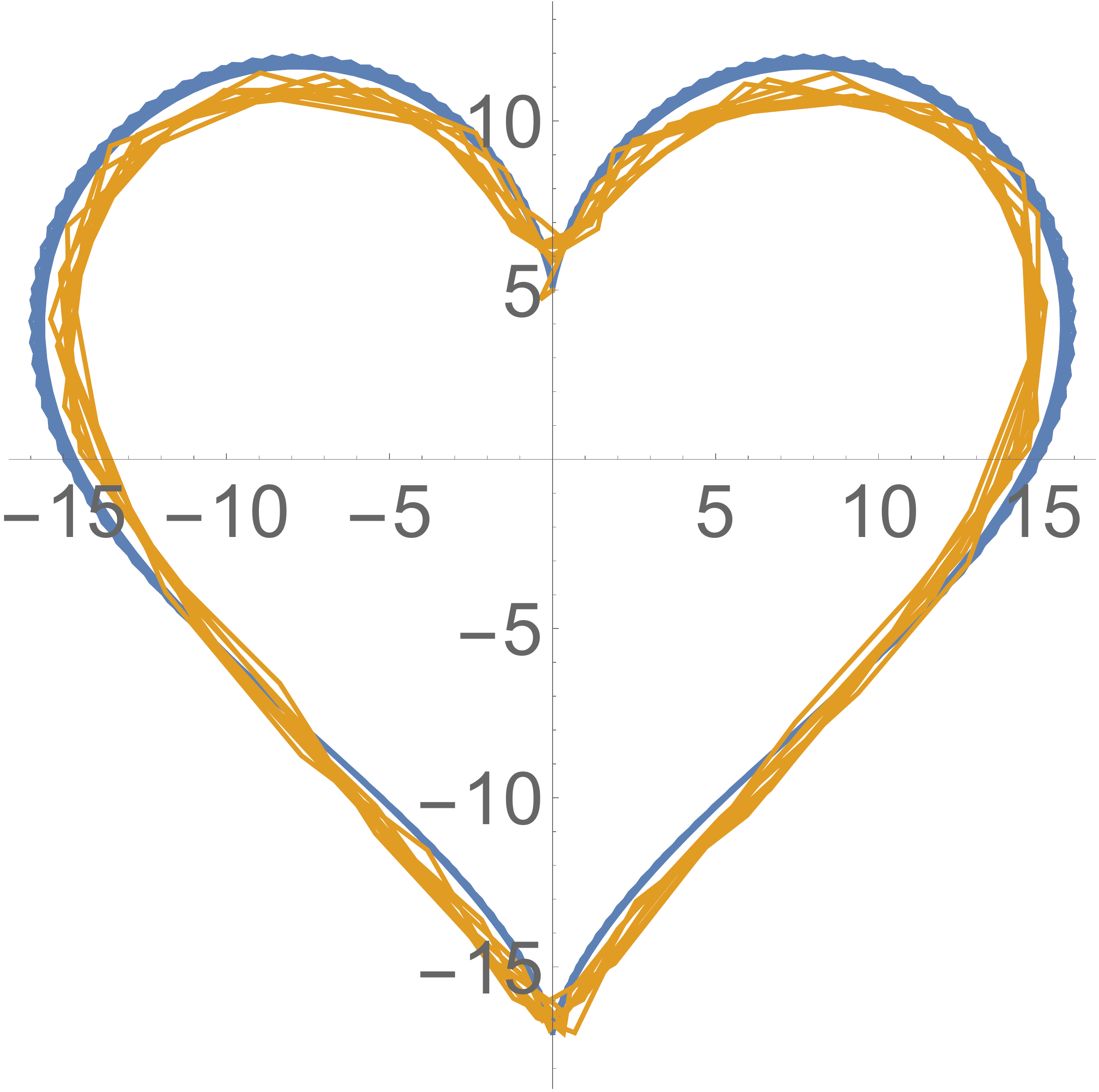}}
  \subcaptionbox{Relative cost: $\Alg / \Opt - 1$.}{\includegraphics[width = 0.25 \hsize]{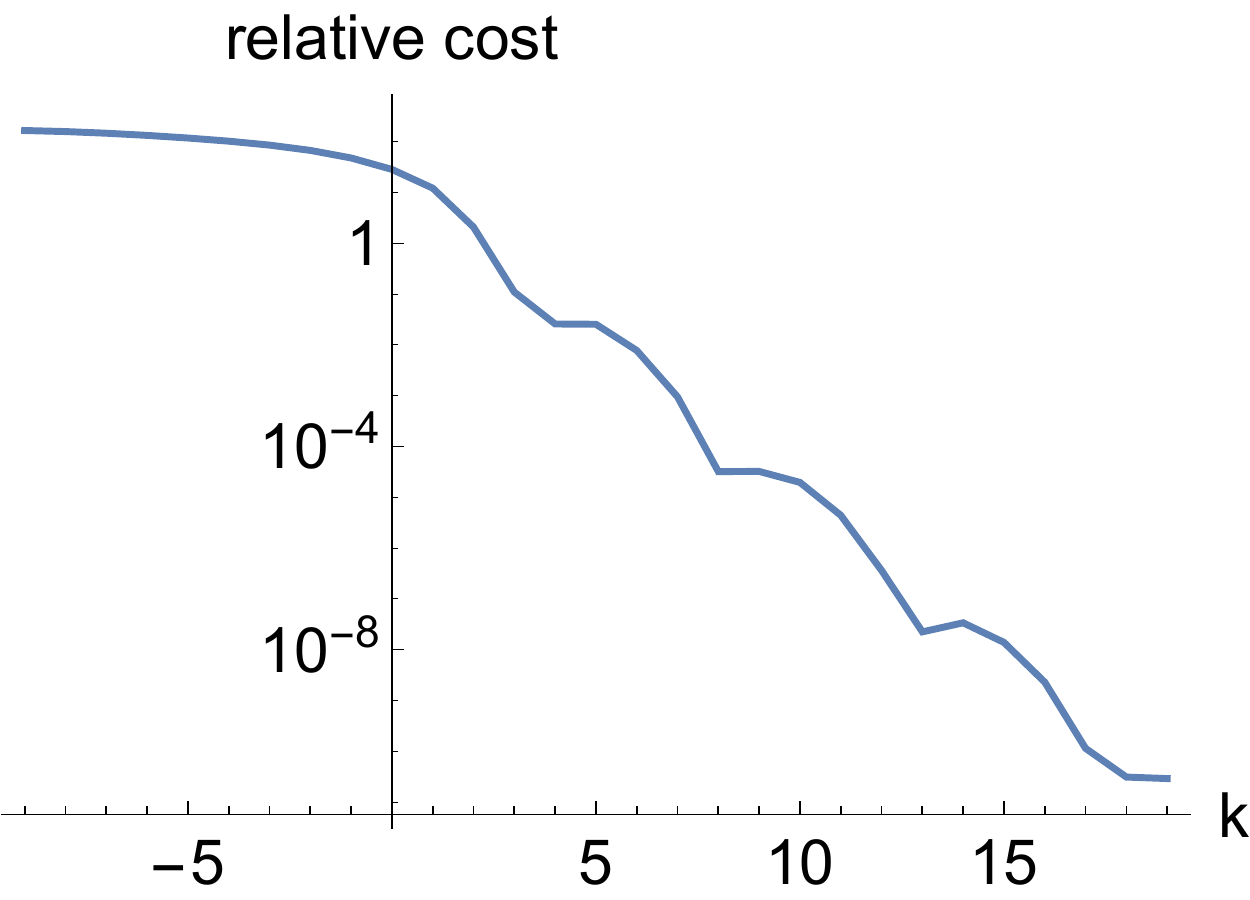}}
  \caption{Tracking results with $k$ exact predictions and no delay, or with $d$ steps of delay with no predictions. (a-d) show the desired trajectory (blue) and actual (orange) trajectories. (e) shows both delay and prediction on the x-axis, with the negative part corresponding to delay.  The y-axis is in log-scale.}
  \label{fig:tracking_results}
\end{figure*}

To illustrate our results, we end the paper with numerical examples that highlight the impact of delayed, inexact predictions. To that end, we consider a 2-d tracking problem with the following trajectory \cite{li2019online}, illustrated in Figure \ref{fig:tracking_results}:
\[y_t = \begin{bmatrix} 16 \sin ^3 (\frac{t}{4}) \\ 13 \cos (\frac{t}{4}) - 5 \cos (\frac{2 t}{4}) - 2 \cos (\frac{3 t}{4}) - \cos (\frac{4 t}{4}) \end{bmatrix}.\]
We consider following double integrator dynamics:
\[p_{t+1} = p_t + 0.2 v_t + h_t, \quad v_{t+1} = v_t + 0.2 u_t + \eta_t,\]
where $p_t \in \R^2$ is the position, $v_t$ is the velocity, $u_t$ is the control, and $h_t, \eta_t \sim \Unif[-1, 1]^2$ are i.i.d.\ noises.
The objective is to minimize
$\sum_{t=0}^{T-1} \norm{p_t - y_t}^2 + 0.0016 \norm{u_t}^2$,
where we let $T = 200$.
This problem can be converted to the standard LQR with disturbance $w_t$ by letting $x_t = \begin{bsmallmatrix} p_t \\ v_t \end{bsmallmatrix}$ and $\tilde w_t = \begin{bsmallmatrix} h_t \\ \eta_t \end{bsmallmatrix}$ and then using the reduction in the LQ tracking example in \Cref{section:model}. 
Note that the disturbances are the combination of a deterministic trajectory and i.i.d.\ noise. In contrast, our theoretical results focus on more challenging adversarial disturbances. Nonetheless, the numerical results are consistent with our theorems.

In our first experiment, we study the effect of the number of delays or predictions. For simplicity, we exclude the effect of inexactness of the predictions --- a prediction is either exact ($\epsilon_i = 0$) or uninformative ($\hat w_{t-d+i|t} = 0$). In this case, each exact prediction cancels a step of delay so, delays can be viewed as ``negative'' predictions.
\Cref{fig:tracking_results} shows the performance of the proposed myopic policy in \Cref{sec:myopic_policy} with different numbers of predictions or delays. We see that the cost exponentially decreases (increases) as the number of predictions (delays) increases.

\begin{figure}[t]
  \centering
  \includegraphics[width = 0.9 \hsize]{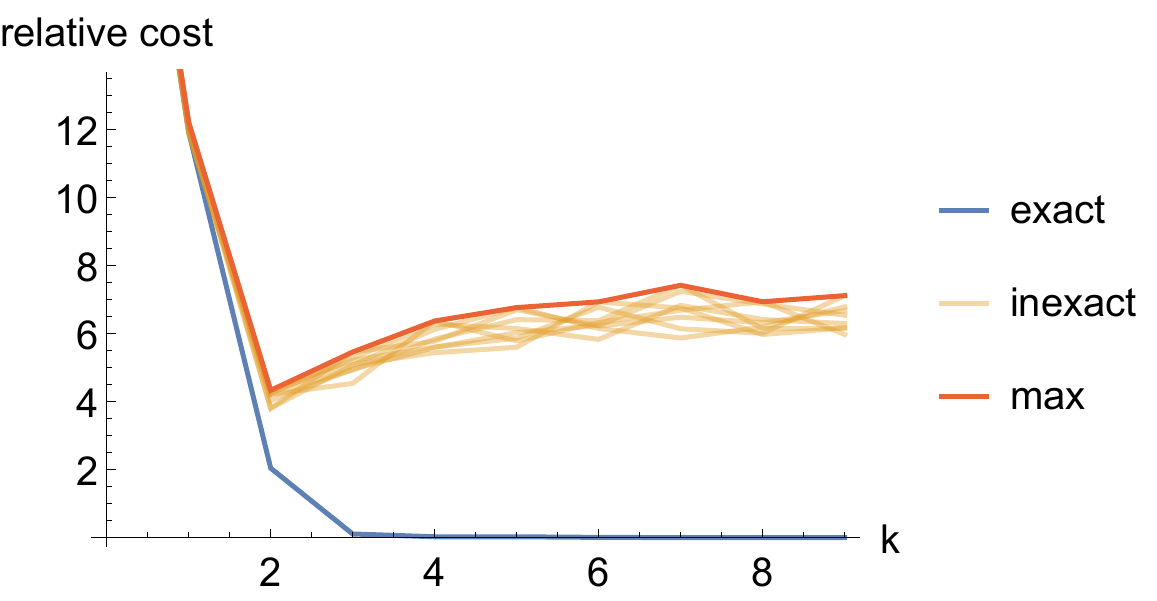}
  \caption{The impact of inexact predictions. The relative cost ($\Alg / \Opt - 1$) of MPC using $k$ {exact} (blue) or {inexact} (orange) predictions is shown.
  }
  \label{fig:inexact}
\end{figure}

In the second experiment, we study the effect of inexact predictions and show that the controller needs to optimize how many predictions are used --- it is better to use only a few predictions and ignore those that are too noisy.  Specifically, we let $\epsilon_i = \frac{1}{5} i^2$, i.e., the noise level of predictions grows quadratically fast with the number of steps into the future. Each estimation error $e_{t+i|t} = w_{t+i} - \hat w_{t+i|t}$ is independently sampled from $\Unif[- \epsilon_i \norm{w_{t+i}}, \epsilon_i \norm{w_{t+i}}]^4$. This process is repeated 8 times, with each instance depicted by an orange line and their maximum represented by a red line. As shown in \Cref{fig:inexact}, with exact predictions, the cost will decrease as the number of predictions increase (the blue line); while with inexact predictions, using fewer predictions may yield better performance.

\section{CONCLUDING REMARKS}
Our results presens the first constant-competitive policy for general LQR control with adversarial disturbances and delayed imperfect information. We also show that in the case of exact predictions with no delay, or in the case of delay with no predictions, the competitive ratio bounds of the proposed myopic policy match the lower bound. However, in the inexact prediction case, the tightness of $\epsilon_i$ in our bounds remains as an open question.
Other important extensions include nonlinear dynamics and time-variant linear systems, which can also lead to studying online learning of robust controllers under model mismatch.

\bibliographystyle{IEEEtran}
\bibliography{IEEEabrv,ref}

\onecolumn
\appendix

\subsection{Cost Characterization Lemma} \label{section:cost_char}

Before we start our proofs, we first present a technical lemma that is used in many of the proofs below. This lemma connects the control cost of a policy to its difference from the offline optimal policy.
\begin{lemma}[Cost characterization] \label{lemma:cost_char_Qf=P}
  Suppose at each time $t$, the controller applies the following policy:
  \begin{equation} \label{eq:affine_policy}
    u_t = -(R + B^\trp P B)^{-1} B^\trp \left(P A x_t + \sum_{i=0}^{T-t-1} {F^\trp}^i P w_{t+i} - \eta_t\right).
  \end{equation}
  If $\Qf = P$, then the control cost is given by:
  \begin{multline} \label{eq:cost_when_Qf=P}
    \Alg = \sum_{t=0}^{T-1} \Bigg(w_t^\trp P w_t + 2 w_t^\trp \sum_{i=1}^{T-t-1} {F^\trp}^i P w_{t+i} - \left(\sum_{i=0}^{T-t-1} {F^\trp}^i P w_{t+i} \right)^\trp H \left( \sum_{i=0}^{T-t-1} {F^\trp}^i P w_{t+i} \right) \Bigg) + \sum_{t=0}^{T-1} \eta_t^\trp H \eta_t \\
    + x_0^\trp P x_0 + 2 x_0^\trp \sum_{i=0}^{T-1} {F^\trp}^{i+1} P w_i.
  \end{multline}
\end{lemma}

Note that the optimal offline policy has $\eta_t = 0$ for all $t$ in \eqref{eq:cost_when_Qf=P} (as derived by \cite{goel2020power}), and as a result, the extra cost of $\Alg$ is given by
\begin{equation} \label{eq:Alg-Opt=eta}
  \Alg - \Opt = \sum_{t=0}^{T-1} \eta_t^\trp H \eta_t.
\end{equation}
In \eqref{eq:affine_policy}, $\eta_t$ can be regarded as the difference between the applied policy and the offline optimal policy. 

We also present below a lemma that has appeared in the body, as we will prove the two lemmas at one time.

\textit{\Cref{lemma:cost_gap_O(1)}.}
\textit{For any algorithm, $\Alg(\Qf) - \Opt^P(\Qf) \le \Alg(P) - \Opt^P(P) + O(1)$, where the $O(1)$ term is with respect to $T$ and it is zero when $\Qf = P$.}

\noindent\hspace{2em}{\itshape Proof of \Cref{lemma:cost_char_Qf=P,lemma:cost_gap_O(1)}:}
Given a disturbance sequence $\vw$, we define the cost-to-go function of a policy described by \eqref{eq:affine_policy}:
\[V_t^\Alg(x_t; \vw)
  \coloneqq \sum_{i=t}^{T-1} (x_i^\trp Q x_i + u_i^\trp R u_i) + x_T^\trp P_T x_T
  = x_t^\trp Q x_t + u_t^\trp R u_t + V_{t+1}^\Alg(x_{t+1}; \vw).\]
We will show by backward induction that $V_t^\Alg(x_t; \vw) = x_t^\trp P_t x_t + x_t^\trp v_t + q_t$ for some $P_t$, $v_t$ and $q_t$. Let $\Delta_t = P_t - P$, where $P$ is the solution of DARE \eqref{eq:DARE}. When $t = T$, we have $P_T = \Qf$, $v_T = 0$ and $q_T = 0$.
Assume this is true at $t + 1$. Then,
\begin{align*}
  & V_t^\Alg(x_t; \vw) \\
  & = x_t^\trp Q x_t + u_t^\trp R u_t + (A x_t + B u_t + w_t)^\trp P_{t+1} (A x_t + B u_t + w_t) + (A x_t + B u_t + w_t)^\trp v_{t+1} + q_{t+1} \\
  & = u_t^\trp (R + B^\trp P_{t+1} B) u_t + 2 u_t^\trp B^\trp (P_{t+1} A x_t + P_{t+1} w_t + v_{t+1}/2) \\
  & \quad + x_t^\trp Q x_t + (A x_t + w_t)^\trp P_{t+1} (A x_t + w_t) + (A x_t + w_t)^\trp v_{t+1} + q_{t+1} \\
  & = u_t^\trp (R + B^\trp P B) u_t + u_t^\trp B^\trp \Delta_{t+1} B u_t \\
  & \quad + 2 u_t^\trp B^\trp (P A x_t + P w_t + v_{t+1}/2) + 2 u_t^\trp B^\trp (\Delta_{t+1} A x_t + \Delta_{t+1} w_t + v_{t+1}/2) \\
  & \quad + x_t^\trp Q x_t + (A x_t + w_t)^\trp P (A x_t + w_t) + (A x_t + w_t)^\trp \Delta_{t+1} (A x_t + w_t) + (A x_t + w_t)^\trp v_{t+1} + q_{t+1} \\
  & = \left(P A x_t + \sum_{i=0}^{T-t-1} {F^\trp}^i P w_{t+i} - \eta_t \right)^\trp H \left(P A x_t + \sum_{i=0}^{T-t-1} {F^\trp}^i P w_{t+i} - \eta_t \right) \\
  & \quad + \left(P A x_t + \sum_{i=0}^{T-t-1} {F^\trp}^i P w_{t+i} - \eta_t \right)^\trp H \Delta_{t+1} H \left(P A x_t + \sum_{i=0}^{T-t-1} {F^\trp}^i P w_{t+i} - \eta_t \right) \\   
  & \quad - 2 \left(P A x_t + \sum_{i=0}^{T-t-1} {F^\trp}^i P w_{t+i} - \eta_t \right)^\trp H (P A x_t + P w_t + v_{t+1}/2) \\
  & \quad - 2 \left(P A x_t + \sum_{i=0}^{T-t-1} {F^\trp}^i P w_{t+i} - \eta_t \right)^\trp H (\Delta_{t+1} A x_t + \Delta_{t+1} w_t + v_{t+1}/2) \\
  & \quad + x_t^\trp Q x_t + (A x_t + w_t)^\trp P (A x_t + w_t) + (A x_t + w_t)^\trp \Delta_{t+1} (A x_t + w_t) + (A x_t + w_t)^\trp v_{t+1} + q_{t+1} \\
  & = (P A x_t)^\trp H (P A x_t) + 2 (P A x_t)^\trp H \left(\sum_{i=0}^{T-t-1} {F^\trp}^i P w_{t+i} - \eta_t \right) \\
  & \quad + \left(\sum_{i=0}^{T-t-1} {F^\trp}^i P w_{t+i} - \eta_t \right)^\trp H \left(\sum_{i=0}^{T-t-1} {F^\trp}^i P w_{t+i} - \eta_t \right) \\
  & \quad + (P A x_t)^\trp H \Delta_{t+1} H (P A x_t) + 2 (P A x_t)^\trp H \Delta_{t+1} H \left(\sum_{i=0}^{T-t-1} {F^\trp}^i P w_{t+i} - \eta_t \right) \\
  & \quad + \left(\sum_{i=0}^{T-t-1} {F^\trp}^i P w_{t+i} - \eta_t \right)^\trp H \Delta_{t+1} H \left(\sum_{i=0}^{T-t-1} {F^\trp}^i P w_{t+i} - \eta_t \right) \\
  & \quad - 2 \left(P A x_t\right)^\trp H (P A x_t + P w_t + v_{t+1}/2) - 2 \left(\sum_{i=0}^{T-t-1} {F^\trp}^i P w_{t+i} - \eta_t \right)^\trp H (P A x_t + P w_t + v_{t+1}/2) \\
  & \quad - 2 \left(P A x_t\right)^\trp H \Delta_{t+1} (A x_t + w_t) - 2 \left(\sum_{i=0}^{T-t-1} {F^\trp}^i P w_{t+i} - \eta_t \right)^\trp H \Delta_{t+1} (A x_t + w_t) \\
  & \quad + x_t^\trp Q x_t + (A x_t + w_t)^\trp P (A x_t + w_t) + (A x_t + w_t)^\trp \Delta_{t+1} (A x_t + w_t) + (A x_t + w_t)^\trp v_{t+1} + q_{t+1} \\
  & = x_t^\trp (Q + A^\trp P A - A^\trp P H P A + F^\trp \Delta_{t+1} F) x_t \\
  & \quad + 2 x_t^\trp F^\trp P_{t+1} w_t + x_t^\trp F^\trp v_{t+1} -  x_t^\trp F^\trp \Delta_{t+1} H \left( \sum_{i=0}^{T-t-1} {F^\trp}^i P w_{t+i} - \eta_t \right) \\
  & \quad + \left(\sum_{i=0}^{T-t-1} {F^\trp}^i P w_{t+i} - \eta_t \right)^\trp H \left(\sum_{i=0}^{T-t-1} {F^\trp}^i P w_{t+i} - \eta_t \right) - 2 \left(\sum_{i=0}^{T-t-1} {F^\trp}^i P w_{t+i} - \eta_t \right)^\trp H (P w_t + v_{t+1}/2) \\
  & \quad + w_t^\trp P w_t + w_t^\trp v_{t+1} + q_{t+1} + O(\Delta_{t+1}).
\end{align*}
Thus, $P + \Delta_t = P_t = Q + A^\trp P A - A^\trp P H P A + F^\trp \Delta_{t+1} F = P + F^\trp \Delta_{t+1} F$ and thus $\Delta_t = F^\trp \Delta_{t+1} F = O(\lambda^{2(T-t)})$, where $\lambda = \frac{1 + \rho(F)}{2}$. The recursive formulae for $v_t$ and $q_t$ are given by:
\begin{align*}
  v_t & = 2 F^\trp P w_t + F^\trp v_{t+1} + O(\lambda^{2(T-t)}) = 2 \sum_{i=0}^{T-t-1} {F^\trp}^{i+1} P w_{t+i} + O(\lambda^{T-t}) \\
  v_{t+1} & = 2 \sum_{i=1}^{T-t-1} {F^\trp}^i P w_{t+i} + O(\lambda^{T-t}), \\
  q_t & = q_{t+1} + w_t^\trp P w_t + 2 w_t^\trp \sum_{i=1}^{T-t-1} {F^\trp}^i P w_{t+i} + \left(\sum_{i=0}^{T-t-1} {F^\trp}^i P w_{t+i} - \eta_t \right)^\trp H \left(\sum_{i=0}^{T-t-1} {F^\trp}^i P w_{t+i} - \eta_t \right) \\
  & \quad - 2 \left(\sum_{i=0}^{T-t-1} {F^\trp}^i P w_{t+i} - \eta_t \right)^\trp H \left(\sum_{i=0}^{T-t-1} {F^\trp}^i P w_{t+i} + O(\lambda^{T-t}) \right) + O(\lambda^{T-t}) \\
  & = q_{t+1} + w_t^\trp P w_t + 2 w_t^\trp \sum_{i=1}^{T-t-1} {F^\trp}^i P w_{t+i} - \left(\sum_{i=0}^{T-t-1} {F^\trp}^i P w_{t+i} \right)^\trp H \left(\sum_{i=0}^{T-t-1} {F^\trp}^i P w_{t+i}\right) + \eta_t^\trp H \eta_t \\
  & \quad + O(\lambda^{T-t}).
\end{align*}
Then,
\[\begin{split}
  \Alg & = V_0^\Alg(x_0; \vw) = x_0^\trp P_0 x_0 + x_0^\trp v_0 + q_0 \\
  & = x_0^\trp P_0 x_0 + x_0^\trp \left( 2 \sum_{i=0}^{T-1} {F^\trp}^{i+1} P w_i + O(\lambda^T) \right) + \sum_{t=0}^{T-1} \eta_t^\trp H \eta_t + \sum_{t=0}^{T-1} O(\lambda^{T-t}) \\
  & \quad + \sum_{t=0}^{T-1} \Bigg(w_t^\trp P w_t + 2 w_t^\trp \sum_{i=1}^{T-t-1} {F^\trp}^i P w_{t+i} - \left( \sum_{i=0}^{T-t-1} {F^\trp}^i P w_{t+i} \right)^\trp H \left(\sum_{i=0}^{T-t-1} {F^\trp}^i P w_{t+i}\right)\Bigg).
\end{split}\]
If $\Qf = P$, then the above $O(\lambda^T)$ and $O(\lambda^{T-t})$ are both zero and thus we obtain \eqref{eq:cost_when_Qf=P}. Otherwise,
\[\Alg(\Qf) - \Opt^P(\Qf) = x_0^\trp O(\lambda^T) + \sum_{t=0}^{T-1} \eta_t^\trp H \eta_t + O(1).\]
Therefore, $(\Alg(\Qf) - \Opt^P(\Qf)) - (\Alg(P) - \Opt^P(P)) = O(1)$. %
\endproof

\subsection{Proof of Theorem \ref{thm:myopic_policy_delay}}

\textit{\Cref{thm:myopic_policy_delay}.}
\textit{Suppose there are $d$ delays and $k$ exact predictions with $k < d$. Assume all used predictions are exact and other disturbances (with unused predictions) are zero. The optimal policy at time $t$ is:
\begin{equation} \tag{\ref{eq:myopic_policy_k_lt_d}}
  u_t = -(R + B^\trp P B)^{-1} B^\trp P A \left(A^{d-k} \hat x_{t-d+k|t} + \sum_{i=0}^{d-k-1} A^i B u_{t-1-i}\right).
\end{equation}}

\begin{proof}
\Cref{lemma:cost_char_Qf=P} implies that when $\Qf = P$, the offline optimal policy is given by
\[u_t = -(R + B^\trp P B)^{-1} B^\trp \left(P A x_t + \sum_{i=0}^{T-t-1} {F^\trp}^i P w_{t+i} \right).\]
However, we are looking for the optimal policy using the incorrect assumptions that (i) $w_{t-d+k}$ and all later disturbances are zero, and (ii) $w_{t-d}, \dots, w_{t-d+k-1}$ equals to $\hat w_{t-d|t}, \dots, \hat w_{t-d+k-1|t}$ respectively.

Replacing $w_{t+i}$ by zero and $x_t$ by $\hat x_{t|t}$ in the above policy, we obtain:
\begin{equation} \label{eq:myopic_policy_delay_implicit}
  u_t = -(R + B^\trp P B)^{-1} B^\trp P A \hat x_{t|t}.
\end{equation}
\[\begin{split}
  \hat x_{t|t} & = A \hat x_{t-1|t} + B u_{t-1} \\
  & = A (A \hat x_{t-2|t} + B u_{t-2}) + B u_{t-1} \\
  & \quad \vdots \\
  & = A^{d-k} \hat x_{t-d+k|t} + A^{d-k-1} B u_{t-d+k} + \dots + B u_{t-1} \\
  & = A^{d-k} \hat x_{t-d+k|t} + \sum_{i=0}^{d-k-1} A^i B u_{t-1-i}.
\end{split}\]
As such, we obtain \Cref{thm:myopic_policy_delay}.
\end{proof}

\subsection{Proof of Theorem \ref{thm:delay&inexact_pred}}

\textit{\Cref{thm:delay&inexact_pred}.}
\textit{Let $c = \norm{P} \|P^{-1}\| (1 + \norm{F})$ and $H = B (R + B^\trp P B)^{-1} B^\trp$. Suppose there are $d$ steps of delays and the controller uses $k$ predictions. When $k \ge d$,
\[\Alg \le \Bigg[\frac{\left(c \sum\limits_{i=0}^{d-1} \epsilon_i \|A^{d-i}\| + c \sum\limits_{i=d}^{k-1} \epsilon_i \|F^{i-d}\| + \|F^{k-d}\| \right)^2}{\norm{H}^{-1} \lambda_{\min}(P^{-1} - F P^{-1} F^\trp - H)} + 1 \Bigg] \Opt + O(1).\]
When $k \le d$,
\[\Alg \le \Bigg[\frac{\bigg(c \sum\limits_{i=0}^{k-1} \epsilon_i \|A^{d-i}\| + c \sum\limits_{i=k}^{d-1} \|A^{d-i}\| + 1 \bigg)^2}{\norm{H}^{-1} \lambda_{\min}(P^{-1} - F P^{-1} F^\trp - H)} + 1 \Bigg] \Opt + O(1).\]
The $O(1)$ is with respect to $T$. It may depend on the system parameters $A$, $B$, $Q$, $R$, $\Qf$ and the range of disturbances $r$, but not on $T$. When $\Qf = P$ the $O(1)$ is zero.}

The proof outline provided in the body lays out a set of lemmas that, together, prove Theorem \ref{thm:delay&inexact_pred}.  Here, we provide proofs for each of them.

\subsection{Proof of Lemma \ref{lemma:when_Qf=P}}

\textit{\Cref{lemma:when_Qf=P}.}
\textit{Suppose $\Qf = P$.  Then, the conclusion of \Cref{thm:delay&inexact_pred} holds.}

\begin{proof}
This lemma considers the case of $\Qf = P$. \Cref{lemma:cost_char_Qf=P} implies that when $\Qf = P$, the cost of the offline optimal policy is:
\begin{multline*}
  \Opt = \sum_{t=0}^{T-1} \Bigg(w_t^\trp P w_t + 2 w_t^\trp \sum_{i=1}^{T-t-1} {F^\trp}^i P w_{t+i} - \left(\sum_{i=0}^{T-t-1} {F^\trp}^i P w_{t+i} \right)^\trp H \left( \sum_{i=0}^{T-t-1} {F^\trp}^i P w_{t+i} \right) \Bigg) \\
  + x_0^\trp P x_0 + 2 x_0^\trp \sum_{i=0}^{T-1} {F^\trp}^{i+1} P w_i.
\end{multline*}
We consider the following substitution:
\begin{equation} \label{eq:substitution}
  \psi_t = \sum_{i=0}^{T-t-1} {F^\trp}^i P w_{t+i}, \quad w_t = P^{-1} (\psi_t - F^\trp \psi_{t+1}).
\end{equation}
Then, the offline optimal cost can be lower bounded:
\begin{equation} \label{eq:Opt_bound}
  \begin{split}
    \Opt & = \sum_{t=0}^{T-1} (w_t^\trp P w_t + 2 w_t^\trp F^\trp \psi_{t+1} - \psi_t^\trp H \psi_t) + x_0^\trp P x_0 + 2 x_0^\trp F^\trp \psi_0 \\
    & = \sum_{t=0}^{T-1} (\psi_t^\trp P^{-1} \psi_t - \psi_{t+1}^\trp F P^{-1} F^\trp \psi_{t+1} - \psi_t^\trp H \psi_t) + x_0^\trp P x_0 + 2 x_0^\trp F^\trp \psi_0 \\
    & = \sum_{t=0}^{T-1} (\psi_t^\trp P^{-1} \psi_t - \psi_t^\trp F P^{-1} F^\trp \psi_t - \psi_t^\trp H \psi_t) + \psi_0^\trp F P^{-1} F^\trp \psi_0 + x_0^\trp P x_0 + 2 x_0^\trp F^\trp \psi_0 \\
    & = \sum_{t=0}^{T-1} \psi_t^\trp (P^{-1} - F P^{-1} F^\trp - H) \psi_t + (F^\trp \psi_0 + P x_0)^\trp P^{-1} (F^\trp \psi_0 + P x_0) \\
    & \ge \lambda_{\min} (P^{-1} - F P^{-1} F^\trp - H) \sum_{t=0}^{T-1} \norm{\psi_t}^2.
  \end{split}
\end{equation}

The myopic policy has two cases and we analyze each of them below.

\paragraph{Case 1: $k \ge d$}
In this case, the controller estimates $x_t$ using $x_{t-d}$ and $\hat w_{t-d|t}, \dots, \hat w_{t-1|t}$.
\[x_t - \hat x_{t|t} = (A x_{t-1} + B u_{t-1} + w_{t-1}) - (A \hat x_{t-1|t} + B u_{t-1} + \hat w_{t-1|t}) = A (x_{t-1} - \hat x_{t-1|t}) + e_{t-1|t}.\]
Applying similar procedures repetitively, we obtain:
\[x_t - \hat x_{t|t} = e_{t-1|t} + A e_{t-2|t} + \dots + A^{d-1} e_{t-d|t} = \sum_{i=1}^d A^{i-1} e_{t-i|t}.\]
Comparing \Cref{eq:myopic_policy_k_gt_d,eq:affine_policy}, we have
\begin{equation} \label{eq:eta_k_gt_d}
  \eta_t = \sum_{i=1}^d P A^i e_{t-i|t} + \sum_{i=0}^{k-d-1} {F^\trp}^i P e_{t+i|t} + \sum_{i=k-d}^{T-t-1} {F^\trp}^i P w_{t+i}.
\end{equation}
Using the substitution in \eqref{eq:substitution}, we bound \eqref{eq:eta_k_gt_d} as follows.
\begin{equation} \label{eq:bound_eta_by_psi_k_gt_d}
\begin{split}
  \norm{\eta_t} & = \norm{\sum_{i=1}^d P A^i e_{t-i|t} + \sum_{i=0}^{k-d-1} {F^\trp}^i P e_{t+i|t} + {F^\trp}^{k-d} \psi_{t+k-d}} \\
  & \le \sum_{i=1}^d \norm{P} \norm{A^i} \epsilon_{d-i} \norm{w_{t-i}} + \sum_{i=0}^{k-d-1} \norm{F^i} \norm{P} \epsilon_{i-d} \norm{w_{t+i}} + \norm{F^{k-d}} \norm{\psi_{t+k-d}} \\
  & \le \sum_{i=1}^d \norm{P} \norm{A^i} \epsilon_{d-i} \norm{P^{-1}} (\norm{\psi_{t-i}} + \norm{F} \norm{\psi_{t-i+1}}) \\
  & \quad + \sum_{i=0}^{k-d-1} \norm{F^i} \norm{P} \epsilon_{i-d} \norm{P^{-1}} (\norm{\psi_{t+i}} + \norm{F} \norm{\psi_{t+i+1}}) + \norm{F^{k-d}} \norm{\psi_{t+k-d}}.
\end{split}
\end{equation}
Note that when $t < d$, some terms in \eqref{eq:bound_eta_by_psi_k_gt_d} have negative subscripts. Those terms do not actually exist and should be regarded as zero. However, for the clarity of the proof, we keep them in the formula. In the later derivations, although we treat them as potentially non-zero, they do not affect our result because we are looking for an upper bound.
Let $\bm\eta = (\|\eta_0\|, \dots, \|\eta_{T-1}\|) \in \R^T$ and $\bm\psi = (\|\psi_0\|, \dots, \|\psi_{T-1}\|)$. \Cref{eq:bound_eta_by_psi_k_gt_d} provides a linear inequality relationship between $\bm\eta$ and $\bm\psi$. We define matrix $M = \{M_{t, s}\}_{t, s = 0}^{T-1} \in \R^{T \times T}$ such that $M_{t, s}$ is the coefficient of $\|\psi_s\|$ in the bound of $\|\eta_t\|$ in \eqref{eq:bound_eta_by_psi_k_gt_d}. Then, $\bm\eta \le M \bm\psi$.
\begin{equation} \label{eq:the_eta_term}
  \sum_{t=0}^{T-1} \eta_t^\trp H \eta_t \le \norm{H} \bm\eta^\trp \bm\eta \le \norm{H} \bm\psi^\trp M^\trp M \bm\psi \le \lambda_{\max}(M^\trp M) \norm{H} \norm{\bm\psi}^2.
\end{equation}
\begin{proposition}[Gershgorin circle theorem] \label{thm:gershgorin}
  Let $A \in \C^{n \times n}$.
  Let $D(A_{i, i}, R_i) \subseteq \C$ be a closed disc centered at $A_{i, i}$ with radius $R_i = \sum_{j \ne i} \abs{A_{i, j}}$. Then, every eigenvalue of $A$ lies within at least one of the discs $D(A_{i, i}, R_i)$.
\end{proposition}
We use \Cref{thm:gershgorin} to bound the eigenvalues of $M^\trp M$:
\begin{equation} \label{eq:lambda_max}
  \lambda_{\max} (M^\trp M) \le \max_i \sum_{j=0}^{T-1} (M^\trp M)_{i, j} = \norm{M^\trp M \bm 1}_\infty.
\end{equation}
Plugging $\norm{\psi_s} = 1$ for all $s$ into \eqref{eq:bound_eta_by_psi_k_gt_d}, we have:
\[M \bm 1 \le \left(\norm{P} \norm{P^{-1}} (1 + \norm{F}) \left(\sum_{i=1}^d \norm{A^i} \epsilon_{d-i} + \sum_{i=0}^{k-d-1} \norm{F^i} \epsilon_{i-d}\right) + \norm{F^{k-d}} \right) \bm 1.\]
Thus, \eqref{eq:lambda_max} can be further bounded:
\[\lambda_{\max}(M^\trp M) \le \left(c \left( \sum_{i=1}^d \norm{A^i} \epsilon_{d-i} + \sum_{i=0}^{k-d-1} \norm{F^i} \epsilon_{i-d} \right) + \norm{F^{k-d}} \right)^2,\]
where $c = \norm{P} \norm{P^{-1}} (1 + \norm{F})$. Together with \Cref{eq:Alg-Opt=eta,eq:the_eta_term}, this implies that
\[\Alg - \Opt = \sum_{t=0}^{T-1} \eta_t^\trp H \eta_t \le \left(c \left( \sum_{i=1}^d \norm{A^i} \epsilon_{d-i} + \sum_{i=0}^{k-d-1} \norm{F^i} \epsilon_{i-d} \right) + \norm{F^{k-d}} \right)^2 \norm{H} \norm{\bm\psi}^2.\]
Together with \eqref{eq:Opt_bound}, we have
\[\frac{\Alg - \Opt}{\Opt} \le \left(c \left( \sum_{i=1}^d \norm{A^i} \epsilon_{d-i} + \sum_{i=0}^{k-d-1} \norm{F^i} \epsilon_{i-d} \right) + \norm{F^{k-d}} \right)^2 \norm{H} \lambda_{\min}^{-1} (P^{-1} - F P^{-1} F^\trp - H).\]

\paragraph{Case 2: $k < d$}

We start from \eqref{eq:myopic_policy_delay_implicit}. In this case, we have the following equations.
\begin{align*}
  x_t - \hat x_{t|t} & = A (x_{t-1} - \hat x_{t-1|t}) + w_{t-1} - 0. \\
  & \vdots \\
  x_{t-d+k+1} - \hat x_{t-d+k+1|t} & = A (x_{t-d+k} - \hat x_{t-d+k|t}) + w_{t-d+k} - 0. \\
  x_{t-d+k} - \hat x_{t-d+k|t} & = A (x_{t-d+k-1} - \hat x_{t-d+k-1|t}) + w_{t-d+k-1} - \hat w_{t-d+k-1|t}. \\
  & \vdots \\
  x_{t-d+1} - \hat x_{t-d+1|t} & = A (x_{t-d} - \hat x_{t-d|t}) + w_{t-d} - \hat w_{t-d|t}.
\end{align*}
Note that in the last line, $x_{t-d} = \hat x_{t-d|t}$. Thus, all of the above equations can be combined into the following:
\[x_t - \hat x_{t|t} = \sum_{i=0}^{k-1} A^{d-i-1} e_{t-d+i|t} + \sum_{i=k}^{d-1} A^{d-i-1} w_{t-d+i}.\]
Therefore, the policy can be written as:
\[\begin{split}
  u_t & = -(R + B^\trp P B)^{-1} B^\trp P A \hat x_{t|t} \\
  & = -(R + B^\trp P B)^{-1} B^\trp P A \left(x_t - \sum_{i=0}^{k-1} A^{d-i-1} e_{t-d+i|t} - \sum_{i=k}^{d-1} A^{d-i-1} w_{t-d+i} \right).
\end{split}\]
We compare this with \eqref{eq:affine_policy} to get
\[\eta_t = \sum_{i=0}^{k-1} P A^{d-i} e_{t-d+i|t} + \sum_{i=k}^{d-1} P A^{d-i} w_{t-d+i} + \sum_{i=0}^{T-t-1} {F^\trp}^i P w_{t+i}.\]
With the substitution in \eqref{eq:substitution},
\begin{equation} \label{eq:bound_eta_by_psi_k_lt_d}
  \begin{split}
    \norm{\eta_t} & = \norm{\sum_{i=0}^{k-1} P A^{d-i} e_{t-d+i|t} + \sum_{i=k}^{d-1} P A^{d-i} w_{t-d+i} + \psi_t} \\
    & \le \sum_{i=0}^{k-1} \norm{P} \norm{A^{d-i}} \epsilon_i \norm{w_{t-d+i}} + \sum_{i=k}^{d-1} \norm{P} \norm{A^{d-i}} \norm{w_{t-d+i}} + \norm{\psi_t} \\
    & \le \sum_{i=0}^{k-1} \norm{P} \norm{A^{d-i}} \epsilon_i \norm{P^{-1}} (\norm{\psi_{t-d+i}} + \norm{F} \norm{\psi_{t-d+i+1}}) \\
    & \quad + \sum_{i=k}^{d-1} \norm{P} \norm{A^{d-i}} \norm{P^{-1}} (\norm{\psi_{t-d+i}} + \norm{F} \norm{\psi_{t-d+i+1}}) + \norm{\psi_t}.
  \end{split}
\end{equation}
Similar to the previous case, we define matrix $M = \{M_{t, s}\}_{t, s = 0}^{T-1} \in \R^{T \times T}$ such that $M_{t, s}$ is the coefficient of $\|\psi_s\|$ in the bound of $\|\eta_t\|$ in \eqref{eq:bound_eta_by_psi_k_lt_d}. Then, by \Cref{thm:gershgorin},
\[\lambda_{\max}(M^\trp M) \le \left(c \sum_{i=0}^{k-1} \norm{A^{d-i}} \epsilon_i + c \sum_{i=k}^{d-1} \norm{A^{d-i}} + \norm{\psi_t} \right)^2.\]
\[\Alg - \Opt = \sum_{t=0}^{T-1} \eta_t^\trp H \eta_t \le \norm{H} \bm\eta^\trp \bm\eta \le \norm{H} \bm\psi^\trp M^\trp M \bm\psi \le \lambda_{\max}(M^\trp M) \norm{H} \norm{\bm\psi}^2.\]
\[\frac{\Alg - \Opt}{\Opt} \le \left(c \sum_{i=0}^{k-1} \norm{A^{d-i}} \epsilon_i + c \sum_{i=k}^{d-1} \norm{A^{d-i}} + \norm{\psi_t} \right)^2 \norm{H} \lambda_{\min}^{-1} (P^{-1} - F P^{-1} F^\trp - H).\]

\end{proof}

\subsection{Proof of Lemma \ref{lemma:cost_gap_O(1)}}
See \Cref{section:cost_char}.

\subsection{Proof of Lemma \ref{lemma:terminal_cost_Opt}}
\textit{\Cref{lemma:terminal_cost_Opt}.}
\textit{The followings are equal up to $O(1)$ difference: $\Opt^P(P),\ \Opt^P(\Qf),\ \Opt^\Qf(\Qf)$.}

\begin{proof}
By definition,
\begin{gather*}
  \Opt^P(P) = \min_X \Opt^X(P) \le \Opt^\Qf(P), \\
  \Opt^\Qf(\Qf) = \min_X \Opt^X(\Qf) \le \Opt^P(\Qf).
\end{gather*}
Moreover, for any $X$,
\[\Opt^X(\Qf) - \Opt^X(P) = x_T^\trp (\Qf - P) x_T = O(1),\]
where $x_T$ is the final state obtained by the policy that is optimal assuming the terminal cost is $X$.
Therefore,
\[\Opt^P(P) \le \Opt^\Qf(P) = \Opt^\Qf(\Qf) + O(1) \le \Opt^P(\Qf) + O(1) = \Opt^P(P) + O(1).\]
As a result, $\Opt^P(P), \Opt^P(\Qf), \Opt^\Qf(\Qf), \Opt^\Qf(P)$ are all equal up to a difference of $O(1)$.\end{proof}

\subsection{Tightness of Theorem \ref{thm:pred_CR}}
The lower bound is obtained in the setting where all disturbances $\{w_t\}$ are i.i.d.\ with zero mean.

Suppose there are $k$ exact predictions and no delays. Let $n = m = 1$.
It has been shown \cite[Theorem 3.2]{yu2020power} that the average cost per time step of the optimal online policy is given by
\[P W - \sum_{i=0}^{k-1} F^{2 i} P^2 H W,\]
where $W$ is the variance of the disturbances.
The minimum offline cost is obtained by taking $k \to \infty$. Thus,
\[\Alg - \Opt = \sum_{i=k}^\infty F^{2 i} P^2 H W = \Theta(F^{2 k}).\]
As a result,
\[\frac{\Alg}{\Opt} = 1 + \Theta(F^{2 k}).\]

This one-dimensional example generalizes to higher dimensions by stacking independent one-dimensional systems together, so that all matrices are diagonal.

\subsection{Tightness of Theorem \ref{thm:delay_CR}}
For the case of $d$ steps of delay and no usable predictions, we can derive a lower bound for the cost of the optimal online policy in the setting of i.i.d.\ noise with zero mean.

Let $n = m = 1$ and assume $\rho(A) > 1$.
Using a similar dynamic programming approach, we can get the cost per time step of the optimal online policy facing $d$ steps of delays, given by
\[A^{2 d} P W + \sum_{i=0}^{d-1} A^{2 i} Q W = \Theta(A^{2d}).\]
As a result,
\[\frac{\Alg}{\Opt} = \Theta(A^{2d}).\]
Similar to the previous example, this one-dimensional example generalizes to high dimensions by simply stacking one-dimensional systems together.

\end{document}